\documentclass{lmcspre}
\usepackage{hyperref}

\usepackage{amsmath}
\usepackage{amssymb}
\usepackage{amsthm}

\usepackage{tikz}
\usetikzlibrary{fadings}

\makeatletter
\newcommand{\douwidehat}[2]{%
	\sbox0{$\m@th#1\widehat{\hphantom{#2}}$}%
	\sbox2{$\m@th#1x$}
	\sbox4{$\m@th#1#2$}
	\dimen0=\ht0
	\advance\dimen0 -.8\ht2
	\dimen2=\dp4
	\rlap{%
		\raisebox{\dimexpr\dimen0-\dimen2}{%
			\scalebox{1}[-1]{\box0}%
		}%
	}%
	{#2}%
}
\makeatother

\title{Uniform Computability of PAC Learning}

\author[V.\ Brattka]{Vasco Brattka\lmcsorcid{0000-0003-4664-2183}}[a]
\author[G.\ Chirache]{Guillaume Chirache\lmcsorcid{0009-0006-9703-8399}}[b]
\address{Faculty of Computer Science, Universit\"at der Bundeswehr M\"unchen, 85577 Neubiberg, Germany\\
and Department of Mathematics and Applied Mathematics, University of Cape Town, Rondebosch 7700, South Africa}
\email{Vasco.Brattka@cca-net.de}
\address{\'Ecole polytechnique, 91128 Palaiseau Cedex, France
}
\email{guillaume.chirache.2022@polytechnique.org}

\begin{document}



\def\AA{{\mathcal A}}
\def\BB{{\mathcal B}}
\def\CC{{\mathcal C}}
\def\DD{{\mathcal D}}
\def\EE{{\mathcal E}}
\def\FF{{\mathcal F}}
\def\GG{{\mathcal G}}
\def\HH{{\mathcal H}}
\def\II{{\mathcal I}}
\def\JJ{{\mathcal J}}
\def\KK{{\mathcal K}}
\def\LL{{\mathcal L}}
\def\MM{{\mathcal M}}
\def\NN{{\mathcal N}}
\def\OO{{\mathcal O}}
\def\PP{{\mathcal P}}
\def\QQ{{\mathcal Q}}
\def\RR{{\mathcal R}}
\def\SS{{\mathcal S}}
\def\TT{{\mathcal T}}
\def\UU{{\mathcal U}}
\def\VV{{\mathcal V}}
\def\WW{{\mathcal W}}
\def\XX{{\mathcal X}}
\def\YY{{\mathcal Y}}
\def\ZZ{{\mathcal Z}}


\def\bA{{\mathbf A}}
\def\bB{{\mathbf B}}
\def\bC{{\mathbf C}}
\def\bD{{\mathbf D}}
\def\bE{{\mathbf E}}
\def\bF{{\mathbf F}}
\def\bG{{\mathbf G}}
\def\bH{{\mathbf H}}
\def\bI{{\mathbf I}}
\def\bJ{{\mathbf J}}
\def\bK{{\mathbf K}}
\def\bL{{\mathbf L}}
\def\bM{{\mathbf M}}
\def\bN{{\mathbf N}}
\def\bO{{\mathbf O}}
\def\bP{{\mathbf P}}
\def\bQ{{\mathbf Q}}
\def\bR{{\mathbf R}}
\def\bS{{\mathbf S}}
\def\bT{{\mathbf T}}
\def\bU{{\mathbf U}}
\def\bV{{\mathbf V}}
\def\bW{{\mathbf W}}
\def\bX{{\mathbf X}}
\def\bY{{\mathbf Y}}
\def\bZ{{\mathbf Z}}


\def\IB{{\mathbb{B}}}
\def\IC{{\mathbb{C}}}
\def\IF{{\mathbb{F}}}
\def\IN{{\mathbb{N}}}
\def\IP{{\mathbb{P}}}
\def\IQ{{\mathbb{Q}}}
\def\IR{{\mathbb{R}}}
\def\IS{{\mathbb{S}}}
\def\IT{{\mathbb{T}}}
\def\IZ{{\mathbb{Z}}}

\def\IIB{{\mathbb{\mathbf B}}}
\def\IIC{{\mathbb{\mathbf C}}}
\def\IIN{{\mathbb{\mathbf N}}}
\def\IIQ{{\mathbb{\mathbf Q}}}
\def\IIR{{\mathbb{\mathbf R}}}
\def\IIZ{{\mathbb{\mathbf Z}}}


\def\ELSE{\quad\mbox{else}\quad}
\def\WITH{\quad\mbox{with}\quad}
\def\FOR{\quad\mbox{for}\quad}
\def\AND{\;\mbox{and}\;}
\def\OR{\;\mbox{or}\;}

\def\To{\longrightarrow}
\def\TO{\Longrightarrow}
\def\In{\subseteq}
\def\sm{\setminus}
\def\Inneq{\In_{\!\!\!\!/}}
\def\dmin{\mathop{\dot{-}}}
\def\splus{\oplus}
\def\SEQ{\triangle}
\def\DIV{\uparrow}
\def\INV{\leftrightarrow}
\def\SET{\Diamond}

\def\kto{\equiv\!\equiv\!>}
\def\kin{\subset\!\subset}
\def\pto{\leadsto}
\def\into{\hookrightarrow}
\def\onto{\to\!\!\!\!\!\to}
\def\prefix{\sqsubseteq}
\def\rel{\leftrightarrow}
\def\mto{\rightrightarrows}

\def\B{{\mathsf{{B}}}}
\def\D{{\mathsf{{D}}}}
\def\G{{\mathsf{{G}}}}
\def\E{{\mathsf{{E}}}}
\def\J{{\mathsf{{J}}}}
\def\K{{\mathsf{{K}}}}
\def\L{{\mathsf{{L}}}}
\def\R{{\mathsf{{R}}}}
\def\T{{\mathsf{{T}}}}
\def\U{{\mathsf{{U}}}}
\def\W{{\mathsf{{W}}}}
\def\Z{{\mathsf{{Z}}}}
\def\w{{\mathsf{{w}}}}
\def\HP{{\mathsf{{H}}}}
\def\C{{\mathsf{{C}}}}
\def\Tot{{\mathsf{{Tot}}}}
\def\Fin{{\mathsf{{Fin}}}}
\def\Cof{{\mathsf{{Cof}}}}
\def\Cor{{\mathsf{{Cor}}}}
\def\Equ{{\mathsf{{Equ}}}}
\def\Com{{\mathsf{{Com}}}}
\def\Inf{{\mathsf{{Inf}}}}
\def\bin{{\mathsf{{bin}}}}

\def\Tr{{\mathrm{Tr}}}
\def\Sierp{{\mathrm Sierpi{\'n}ski}}
\def\psisierp{{\psi^{\mbox{\scriptsize\Sierp}}}}
\def\cl{{\mathrm{{cl}}}}
\def\Haus{{\mathrm{{H}}}}
\def\Ls{{\mathrm{{Ls}}}}
\def\Li{{\mathrm{{Li}}}}

\def\CL{\mathsf{CL}}
\def\ACC{\mathsf{ACC}}
\def\DNC{\mathsf{DNC}}
\def\ATR{\mathsf{ATR}}
\def\LPO{\mathsf{LPO}}
\def\LLPO{\mathsf{LLPO}}
\def\WKL{\mathsf{WKL}}
\def\WWKL{\mathsf{WWKL}}
\def\RCA{\mathsf{RCA}}
\def\ACA{\mathsf{ACA}}
\def\SEP{\mathsf{SEP}}
\def\BCT{\mathsf{BCT}}
\def\IVT{\mathsf{IVT}}
\def\IMT{\mathsf{IMT}}
\def\OMT{\mathsf{OMT}}
\def\CGT{\mathsf{CGT}}
\def\UBT{\mathsf{UBT}}
\def\BWT{\mathsf{BWT}}
\def\HBT{\mathsf{HBT}}
\def\BFT{\mathsf{BFT}}
\def\FPT{\mathsf{FPT}}
\def\WAT{\mathsf{WAT}}
\def\LIN{\mathsf{LIN}}
\def\B{\mathsf{B}}
\def\BF{\mathsf{B_\mathsf{F}}}
\def\BI{\mathsf{B_\mathsf{I}}}
\def\C{\mathsf{C}}
\def\CF{\mathsf{C_\mathsf{F}}}
\def\CN{\mathsf{C_{\IN}}}
\def\CI{\mathsf{C_\mathsf{I}}}
\def\CK{\mathsf{C_\mathsf{K}}}
\def\CA{\mathsf{C_\mathsf{A}}}
\def\WPO{\mathsf{WPO}}
\def\WLPO{\mathsf{WLPO}}
\def\MP{\mathsf{MP}}
\def\BD{\mathsf{BD}}
\def\Fix{\mathsf{Fix}}
\def\Mod{\mathsf{Mod}}

\def\TD{\mathsf{TD}}
\def\Kol{\mathrm{K}}
\def\PA{\mathsf{PA}}
\def\HYP{\mathsf{HYP}}
\def\CFC{\mathsf{CFC}}
\def\BMT{\mathsf{BMT}}
\def\NEQ{\mathsf{NEQ}}
\def\NON{\mathsf{NON}}
\def\NDOM{\mathsf{NDOM}}
\def\NRNG{\mathsf{NRNG}}
\def\SORT{\mathsf{SORT}}
\def\EC{\mathsf{EC}}
\def\DIS{\mathsf{DIS}}
\def\INF{\mathsf{INF}}
\def\EC{\mathsf{EC}}
\def\AC{\mathsf{AC}}
\def\UG{\mathsf{UG}}
\def\AD{\mathsf{AD}}
\def\ZFC{\mathsf{ZFC}}
\def\ZF{\mathsf{ZF}}
\def\DC{\mathsf{DC}}
\def\BP{\mathsf{BP}}
\def\card{\mathrm{card}}
\def\deg{\mathrm{deg}}
\def\r{\mathrm{r}}

\def\PAC{\mathsf{PAC}}
\def\PPAC{\mathsf{PPAC}}
\def\IPAC{\mathsf{IPAC}}
\def\RPAC{\mathsf{RPAC}}
\def\ERM{\mathsf{ERM}}
\def\WIT{\mathsf{WIT}}
\def\VCC{\mathsf{VCC}}

\def\limninftynd{\lim\nolimits_{\IN_\infty}^{\nearrow}}

\def\s{\mathrm{s}}
\def\r{\mathrm{r}}
\def\w{\mathsf{w}}

\def\leqm{\mathop{\leq_{\mathrm{m}}}}
\def\equivm{\mathop{\equiv_{\mathrm{m}}}}
\def\leqT{\mathop{\leq_{\mathrm{T}}}}
\def\lT{\mathop{<_{\mathrm{T}}}}
\def\nleqT{\mathop{\not\leq_{\mathrm{T}}}}
\def\equivT{\mathop{\equiv_{\mathrm{T}}}}
\def\nequivT{\mathop{\not\equiv_{\mathrm{T}}}}
\def\leqwtt{\mathop{\leq_{\mathrm{wtt}}}}
\def\equiPT{\mathop{\equiv_{\P\mathrm{T}}}}
\def\leqW{\mathop{\leq_{\mathrm{W}}}}
\def\equivW{\mathop{\equiv_{\mathrm{W}}}}
\def\leqtW{\mathop{\leq_{\mathrm{tW}}}}
\def\leqSW{\mathop{\leq_{\mathrm{sW}}}}
\def\equivSW{\mathop{\equiv_{\mathrm{sW}}}}
\def\leqPW{\mathop{\leq_{\widehat{\mathrm{W}}}}}
\def\equivPW{\mathop{\equiv_{\widehat{\mathrm{W}}}}}
\def\leqFPW{\mathop{\leq_{\mathrm{W}^*}}}
\def\equivFPW{\mathop{\equiv_{\mathrm{W}^*}}}
\def\leqWW{\mathop{\leq_{\overline{\mathrm{W}}}}}
\def\nleqW{\mathop{\not\leq_{\mathrm{W}}}}
\def\nleqSW{\mathop{\not\leq_{\mathrm{sW}}}}
\def\lW{\mathop{<_{\mathrm{W}}}}
\def\lSW{\mathop{<_{\mathrm{sW}}}}
\def\nW{\mathop{|_{\mathrm{W}}}}
\def\nSW{\mathop{|_{\mathrm{sW}}}}
\def\leqt{\mathop{\leq_{\mathrm{t}}}}
\def\equivt{\mathop{\equiv_{\mathrm{t}}}}
\def\leqtop{\mathop{\leq_\mathrm{t}}}
\def\equivtop{\mathop{\equiv_\mathrm{t}}}

\def\bigtimes{\mathop{\mathsf{X}}}

\def\leqm{\mathop{\leq_{\mathrm{m}}}}
\def\equivm{\mathop{\equiv_{\mathrm{m}}}}
\def\leqT{\mathop{\leq_{\mathrm{T}}}}
\def\leqM{\mathop{\leq_{\mathrm{M}}}}
\def\equivT{\mathop{\equiv_{\mathrm{T}}}}
\def\equiPT{\mathop{\equiv_{\P\mathrm{T}}}}
\def\leqW{\mathop{\leq_{\mathrm{W}}}}
\def\equivW{\mathop{\equiv_{\mathrm{W}}}}
\def\nequivW{\mathop{\not\equiv_{\mathrm{W}}}}
\def\leqSW{\mathop{\leq_{\mathrm{sW}}}}
\def\equivSW{\mathop{\equiv_{\mathrm{sW}}}}
\def\leqPW{\mathop{\leq_{\widehat{\mathrm{W}}}}}
\def\equivPW{\mathop{\equiv_{\widehat{\mathrm{W}}}}}
\def\nleqW{\mathop{\not\leq_{\mathrm{W}}}}
\def\nleqSW{\mathop{\not\leq_{\mathrm{sW}}}}
\def\lW{\mathop{<_{\mathrm{W}}}}
\def\lSW{\mathop{<_{\mathrm{sW}}}}
\def\nW{\mathop{|_{\mathrm{W}}}}
\def\nSW{\mathop{|_{\mathrm{sW}}}}

\def\botW{\mathbf{0}}
\def\midW{\mathbf{1}}
\def\topW{\mathbf{\infty}}

\def\pol{{\leq_{\mathrm{pol}}}}
\def\rem{{\mathop{\mathrm{rm}}}}

\def\cc{{\mathrm{c}}}
\def\d{{\,\mathrm{d}}}
\def\e{{\mathrm{e}}}
\def\ii{{\mathrm{i}}}

\def\Cf{C\!f}
\def\id{{\mathrm{id}}}
\def\pr{{\mathrm{pr}}}
\def\inj{{\mathrm{inj}}}
\def\cf{{\mathrm{cf}}}
\def\dom{{\mathrm{dom}}}
\def\range{{\mathrm{range}}}
\def\graph{{\mathrm{graph}}}
\def\Graph{{\mathrm{Graph}}}
\def\epi{{\mathrm{epi}}}
\def\hypo{{\mathrm{hypo}}}
\def\Lim{{\mathrm{Lim}}}
\def\diam{{\mathrm{diam}}}
\def\dist{{\mathrm{dist}}}
\def\supp{{\mathrm{supp}}}
\def\union{{\mathrm{union}}}
\def\fiber{{\mathrm{fiber}}}
\def\ev{{\mathrm{ev}}}
\def\mod{{\mathrm{mod}}}
\def\sat{{\mathrm{sat}}}
\def\seq{{\mathrm{seq}}}
\def\lev{{\mathrm{lev}}}
\def\mind{{\mathrm{mind}}}
\def\arccot{{\mathrm{arccot}}}
\def\cl{{\mathrm{cl}}}

\def\Add{{\mathrm{Add}}}
\def\Mul{{\mathrm{Mul}}}
\def\SMul{{\mathrm{SMul}}}
\def\Neg{{\mathrm{Neg}}}
\def\Inv{{\mathrm{Inv}}}
\def\Ord{{\mathrm{Ord}}}
\def\Sqrt{{\mathrm{Sqrt}}}
\def\Re{{\mathrm{Re}}}
\def\Im{{\mathrm{Im}}}
\def\Sup{{\mathrm{Sup}}}

\def\LSC{{\mathcal LSC}}
\def\USC{{\mathcal USC}}

\def\CE{{\mathcal{E}}}
\def\Pref{{\mathrm{Pref}}}

\def\Baire{\IN^\IN}

\def\TRUE{{\mathrm{TRUE}}}
\def\FALSE{{\mathrm{FALSE}}}

\def\co{{\mathrm{co}}}

\def\BBB{{\tt B}}

\def\Prob{{\IP}}
\def\VCdim{{\mathrm{VCdim}}}
\def\eVCdim{{\mathrm{eVCdim}}}
\def\argmin{{\mathop{\mathrm{argmin}}}}
\def\divv{{\mathop{\mathrm{div}}}}

\newcommand{\PPACd}[1]{\PPAC_{<#1}}
\newcommand{\IPACd}[1]{\IPAC_{<#1}}
\newcommand{\RPACd}[1]{\RPAC_{<#1}}

\newcommand{\SO}[1]{{{\mathbf\Sigma}^0_{#1}}}
\newcommand{\SI}[1]{{{\mathbf\Sigma}^1_{#1}}}
\newcommand{\PO}[1]{{{\mathbf\Pi}^0_{#1}}}
\newcommand{\PI}[1]{{{\mathbf\Pi}^1_{#1}}}
\newcommand{\DO}[1]{{{\mathbf\Delta}^0_{#1}}}
\newcommand{\DI}[1]{{{\mathbf\Delta}^1_{#1}}}
\newcommand{\sO}[1]{{\Sigma^0_{#1}}}
\newcommand{\sI}[1]{{\Sigma^1_{#1}}}
\newcommand{\pO}[1]{{\Pi^0_{#1}}}
\newcommand{\pI}[1]{{\Pi^1_{#1}}}
\newcommand{\dO}[1]{{{\Delta}^0_{#1}}}
\newcommand{\dI}[1]{{{\Delta}^1_{#1}}}
\newcommand{\sP}[1]{{\Sigma^\P_{#1}}}
\newcommand{\pP}[1]{{\Pi^\P_{#1}}}
\newcommand{\dP}[1]{{{\Delta}^\P_{#1}}}
\newcommand{\sE}[1]{{\Sigma^{-1}_{#1}}}
\newcommand{\pE}[1]{{\Pi^{-1}_{#1}}}
\newcommand{\dE}[1]{{\Delta^{-1}_{#1}}}

\newcommand{\dBar}[1]{{\overline{\overline{#1}}}}

\def\QED{$\hspace*{\fill}\Box$}
\def\rand#1{\marginpar{\rule[-#1 mm]{1mm}{#1mm}}}

\def\BL{\BB}


\newcommand{\bra}[1]{\langle#1|}
\newcommand{\ket}[1]{|#1\rangle}
\newcommand{\braket}[2]{\langle#1|#2\rangle}

\newcommand{\ind}[1]{{\em #1}\index{#1}}
\newcommand{\mathbox}[1]{\[\fbox{\rule[-4mm]{0cm}{1cm}$\quad#1$\quad}\]}


\newenvironment{eqcase}{\left\{\begin{array}{lcl}}{\end{array}\right.}

\theoremstyle{definition}
\newtheorem{theorem}{Theorem}
\newtheorem{definition}[theorem]{Definition}
\newtheorem{problem}[theorem]{Problem}
\newtheorem{assumption}[theorem]{Assumption}
\newtheorem{corollary}[theorem]{Corollary}
\newtheorem{proposition}[theorem]{Proposition}
\newtheorem{lemma}[theorem]{Lemma}
\newtheorem{observation}[theorem]{Observation}
\newtheorem{question}[theorem]{Question}
\newtheorem{example}[theorem]{Example}
\newtheorem{convention}[theorem]{Convention}
\newtheorem{conjecture}[theorem]{Conjecture}

\subjclass{[{\bf Theory of computation}]:  Logic; [{\bf Mathematics of computing}]: Continuous mathematics.}

\newcommand{\gnote}[1]{\textcolor{red}{#1}}

\maketitle              

\begin{abstract}
	We study uniform computability properties of PAC learning using Weihrauch complexity.
	We focus on closed concept classes, which are either represented by positive, by negative
	or by full information.
	Among other results, we prove that proper PAC learning from positive information
	is equivalent to the limit operation on Baire space, whereas improper
	PAC learning from positive information is closely related to Weak K\H{o}nig's Lemma and even equivalent
	to it, when we have some negative information about the admissible hypotheses.
	If arbitrary hypotheses are allowed, then improper PAC learning from positive information
	is still in a finitary DNC range, which implies that it is non-deterministically computable, but does not allow for probabilistic algorithms.
	These results can also be seen as a classification of the degree of constructivity of the
	Fundamental Theorem of Statistical Learning.
	All the aforementioned results hold if an upper bound of the VC dimension is
	provided as an additional input information.
	We also study the question of how these results are affected if the VC dimension
	is not given, but only promised to be finite or if concept classes are represented
	by negative or full information. 
	Finally, we also classify the complexity of the VC dimension operation itself,
	which is a problem that is of independent interest. For positive or full information it 
	turns out to be equivalent to the binary sorting problem, for negative
	information it is equivalent to the jump of sorting.
	This classification allows also conclusions regarding the Borel complexity
	of PAC learnability.\\

\end{abstract}

\section{Introduction}

In recent years computability properties of PAC learning have found increasing
interest in the research community
and several authors have investigated computable variants of PAC
learning, starting from a joint work of Agarwal, Ananthakrishnan, Ben-David, Lechner, and Urner ~\cite{AAB+20,AAD+22,Ste22a,DRKRS23,GTU24,AH25,KK25}.
An early result in this direction was already proved by Soloveichik~\cite{Sol08}.
A good survey on the results can be obtained from the most recent work of
Kattermann and Krapp~\cite{KK25}.
In most of these cases PAC learning is looked at from a non-uniform perspective
and the focus is on computability properties of individual concept classes and their learners.
We take a different standpoint and investigate PAC learning from a uniform
perspective, using the tools of computable analysis~\cite{Wei00,BH21} and
Weihrauch complexity~\cite{BGP21}.
Concepts from computable analysis have already been applied in other work in the 
context of learning theory theory, e.g., in~\cite{dBre10,dBY10a,HK14,Cip23,CMPR23,Bra23a}. Most notably Weihrauch reducibility has already been used to analyze PAC learning
in the continuous setting 
by Ackerman, Asilis, Di, Freer and Tristan~\cite{AAD+22}.

The uniform approach using Weihrauch complexity allows us to provide succinct and precise answers to the
following three questions at once. These questions roughly come from a logical, a computability-theoretical
and a theoretical computer science perspective, respectively:
\begin{enumerate}
	\item How non-constructive is the Fundamental Theorem of Statistical Learning?
	\item What is the computability-theoretic complexity of the problem:
	      {\em given} a concept class of finite VC dimension, {\em find} a learner for it.
	\item Which computational properties does PAC learning have: is it computable,
	      limit computable, non-deterministically computable and so forth?
\end{enumerate}

In order to describe our results
we start recalling some basic definitions of computational learning
theory. For details we refer the reader to standard references such as \cite{SB14,MRT18}.

By $\SS:=(\IN\times\{0,1\})^*$ we denote the {\em set of samples}, which
are finite words of pairs $(x,y)\in\IN\times\{0,1\}$.
A {\em hypothesis} is a function $h:\IN\to\{0,1\}$ and
a {\em learner} is a function $A:\SS\to2^\IN$ that produces a
hypothesis for every sample. We interchangeably view the points in $2^\IN$ either as
binary sequences or as subsets of $\IN$ identified with their characteristic functions.
By $\LL$ we denote the set of all learners.
A {\em hypothesis class} is a set $\HH\In2^\IN$.
Sometimes we work with two hypothesis classes $\CC\In\HH\In2^\IN$,
in which case we call the smaller one $\CC$ also the {\em concept class}.
The goal of PAC learning is to find a learner $A:\SS\to\HH$ for every
suitable concept class $\CC\In\HH$ such that this learner fits well to the concept class.
What it means that the learner fits well to a concept class
is expressed in terms of probability and roughly speaking the goal
is that the learner makes large mistakes only with a small probability
(PAC learning stands for ``probably approximately correct'' learning).

We use {\em probability distributions} $\DD:\IN\times\{0,1\}\to\IR$,
which are simply functions that take non-negative values and satisfy
$\sum_{(n,i)\in\IN\times\{0,1\}}\DD(n,i)=1$.
By $\PP$ we denote the set of such probability distributions.
For a probability distribution $\DD$ and a hypothesis $h$ we denote
the {\em true error (risk)} of $h$ according to $\DD$ by
\[L_\DD(h):=\Prob_{(x,y)\sim\DD}[h(x)\not=y].\]
Here $\Prob_{(x,y)\sim\DD}$ denotes the probability over all
$(x,y)\in\IN\times\{0,1\}$ that
are distributed according to $\DD$.
Likewise, we denote by $\Prob_{S\sim\DD^n}$ the probability
over all samples $S\in\SS$ of length $n$ that are independently
and identically distributed (i.i.d.) according to $\DD$.

Now we are prepared to define what it means that a learner PAC learns
a hypothesis class. We use the {\em Cantor tupling function} with 
$\langle n,k\rangle:=\frac{1}{2}(n+k)(n+k+1)+k$ for all $n,k\in\IN$.

\begin{definition}[PAC learning]
	Let $\varnothing\not=\CC\In\HH\In2^\IN$ be hypothesis classes, let $A:\SS\to\HH$ be a learner
	and let $m:\IN\to\IN$ be a function that we call {\em sample complexity}.
	We say that $A$ {\em PAC learns}
	$\CC$ {\em relative to} $\HH$ {\em with} $m$
	if for all $i,j\in\IN$ and $n\geq m\langle i,j\rangle$
	and for every probability distribution $\DD$ over $\IN\times\{0,1\}$
	the following holds:
	\begin{eqnarray}
		\IP_{S\sim\DD^n}\left[L_\DD(A(S))\leq\inf_{h\in\CC}L_\DD(h)+2^{-i}\right]> 1-2^{-j}.
		\label{eq:PAC}
	\end{eqnarray}
	If $\CC=\HH$, then we say that $A$ {\em (properly) PAC learns} $\CC$ {\em with} $m$. 
	If $\CC\In\HH$ and $\HH=2^\IN$, then we say that $A$ {\em improperly PAC learns} $\CC$ {\em with} $m$.
\end{definition}

The parameter $\varepsilon=2^{-i}$ is sometimes called the {\em accuracy}
and $\delta=2^{-j}$ the {\em confidence}.
Roughly speaking, the condition expresses that over all samples $S$ of sufficient
size $n\geq m\langle i,j\rangle$ with high confidence the true error of the hypothesis
$A(S)\in\HH$ produced by the learner should be  with high accuracy
close to the {\em optimal true error}
$\inf_{h\in\CC} L_{\DD}(h)$ that can be achieved in the concept class $\CC$ (in both
cases ``high'' means large $i$ and $j$, respectively).
Improper PAC learning is also known under the name {\em representation independent learning}~\cite[Remark~3.2]{SB14}. However, we will stick to the shorter and terminology of improper learning.

\begin{definition}[PAC learnability]
	Let $\varnothing\not=\CC\In\HH\In2^\IN$ be hypothesis classes.
	We say that $\CC$ is {\em PAC learnable relative to} $\HH$ if there is a learner
	$A:\SS\to\HH$ and a sample complexity $m:\IN\to\IN$ such that $A$ PAC learns $\CC$ relative to $\HH$ with $m$. 
	If $\CC=\HH$, then we briefly say  that $\CC$ is {\em (properly) PAC learnable}. 
	If $\CC\In\HH$ and $\HH=2^\IN$, then we briefly say that $\CC$ is {\em improperly PAC learnable},
\end{definition}

In cases of doubt PAC learnability always refers to proper PAC learnability.
An obvious question is whether there are any PAC learnable concept classes $\CC$
and, if so, how they can be characterized. This question is fully answered by the
Fundamental Theorem of Statistical Learning, which was proved by
Blumer, Ehrenfeucht, Haussler, and Warmuth~\cite[Theorem~2.1]{BEHW89} (see also~\cite[Theorem~6.7]{SB14}), building on earlier work of
Vapnik and Chervonenkis~\cite{VC71}.
The Fundamental Theorem states that $\CC$ is PAC learnable if and only if its
VC dimension is finite. 

We recall the definition of {\em VC dimension} (also known
as {\em Vapnik Chervonenkis dimension}).
Firstly, we define $\CC\cap A:=\{B\cap A:B\in\CC\}$ for every concept class $\CC\In2^\IN$ and $A\in2^\IN$
and we say that $A$ is {\em shattered} by $\CC$ if $|\CC\cap A|=2^{|A|}$ where $|X|$
denotes the {\em cardinality} of a set $X$.
Using this terminology we can define
\[\VCdim(\CC):=\sup\{|A|\in\IN:A\In\IN\mbox{ is finite and shattered by }\CC\}\in\IN_\infty:=\IN\cup\{\infty\}\]
for every non-empty class $\CC$.
If $\CC=\varnothing$ then the supremum is taken over the empty set and
a natural value would be $-1$. However, in order to stick to the natural numbers, we just
artificially define $\VCdim(\varnothing):=0$.

For classes $\CC$ that are PAC learnable, there is a universal class of PAC learners.
A learner $A:\SS\to\CC$ is called an {\em empirical risk minimizer}
for the class $\CC$ if
\[A(S)\in\argmin_{h\in\CC}L_S(h):=\{h\in\CC:L_S(h)=\min\{L_S(h):h\in\CC\}\}\]
for all $S\in\SS$. Here, we denote by
\[L_S(h):=\frac{|\{i\in\{1,...,n\}:h(x_i)\not=y_i\}|}{n}\]
the so-called {\em empirical risk} of
a sample $S=((x_1,y_1),...,(x_n,y_n))\in\SS$
for a hypothesis $h:\IN\to\{0,1\}$.

We can now summarize the Fundamental Theorem of Statistical Learning in
a version that is relevant for us (see \cite[Theorem~2.1]{BEHW89} and
\cite[Theorems~6.7 and 6.8]{SB14}).

\begin{theorem}[Fundamental Theorem of Statistical Learning]
	\label{thm:fundamental}
	Let $\varnothing\not=\CC\In2^\IN$. Then
	\[\CC\mbox{ is PAC learnable}\iff \VCdim(\CC)\mbox{ is finite.}\]
	In fact, one can make the following additional quantitative statements\footnote{The bounds given in this theorem are typically expressed in terms
	of $\varepsilon=2^{-i}$ and $\delta=2^{-j}$, in which case they read
	\[m_d(\varepsilon,\delta)=c\cdot\frac{d\cdot\log(\frac{1}{\varepsilon})+\log(\frac{1}{\delta})}{\varepsilon}\mbox{ and }M_d(\varepsilon,\delta)=C\cdot\frac{d+\log(\frac{1}{\delta})}{\varepsilon}.\]
	A simple calculation leads to the bounds given in our version of the theorem (with suitably modified values of the constants $c$ and $C$).}
	with fixed constants $c,C>0$ that are independent of $\CC$:
	\begin{enumerate}
		\item If $\VCdim(\CC)\leq d$, then every empirical risk minimizer $A:\SS\to\CC$
		      PAC learns $\CC$ with sample complexity $m$ less or equal than $m_d\langle i,j\rangle:=c\cdot 2^i\cdot(d\cdot i+j)$.
		\item If $d\leq\VCdim(\CC)$, then no learning function $A:\SS\to\HH$ for any $\HH\supseteq\CC$ can
		      PAC learn $\CC$ relative to $\HH$ with sample complexity $m$ less than
		      $M_d\langle i,j\rangle:=C\cdot2^i\cdot(d+j)$.
	\end{enumerate}
\end{theorem}

Without loss of generality, we can assume that $c\in\IN$ and that $C\in\IR$ is computable. This implies that $m_d:\IN\to\IN$ and $M_d:\IN\to\IR$ are computable
functions (also in $d$).

Since we want to analyze the uniform computability-theoretic content of the Fundamental
Theorem of Statistical Learning, we need to inspect its logical form
as a $\forall\exists$--statement. For the most interesting direction ``$\Longleftarrow$''
we essentially obtain the following logical form for the case of proper PAC learning:
\begin{eqnarray*}
(\forall\CC)(\VCdim(\CC)<\infty\TO(\exists (A,m)\in\LL\times\IN^\IN)(\mbox{$A$ PAC learns $\CC$ with $m$}))
\end{eqnarray*}

In order to analyze the Weihrauch complexity of this statement,
we can straightforwardly interpret the corresponding
computational problem as a multivalued map $\PPAC$ that takes $\CC$ as input and produces a possible pair $(A,m)$ as output.
However, this immediately raises a number of questions:

\begin{enumerate}
	\item[(A)] How do we represent $\CC$ and what kind of sets
	      do we want to allow here?
	\item[(B)] Does it make a difference whether an upper bound of $d=\VCdim(\CC)$
	      is provided as an additional input information or whether $\CC$ just comes with the
	      promise that its VC dimension is finite (or has some fixed finite value)?
	\item[(C)] Do we need the sample complexity $m$ on the output side or is $A$
	      also useful without it?
\end{enumerate}

Regarding question (C) it is easy to see that a learner $A$ without corresponding
sample complexity $m$ is pretty useless from a practical perspective. 
We will later briefly touch upon this question again at the end of Section~\ref{sec:IPAC},
where we make one case study,
but we will not systematically study all variants of PAC learning without the
sample complexity as output.

Regarding question (A) we have made the choice to restrict our study to
topologically closed concept classes $\CC$. This choice is reasonable for at least two reasons.
Firstly, the VC dimensions of a set and its closure are identical
and hence closed concept classes are the natural representatives of general
concept classes. 
We do not lose anything essential in terms of PAC learnability by restricting
everything to closed sets. This applies in particular in the improper setting, 
where we do not require that the learner produces hypotheses in the given concept class $\CC$.
Secondly, the class of closed sets has continuum cardinality and can be dealt with using standard techniques of computability theory and computable analysis.
In fact there are three well-established standard ways of representing the space $\AA$ of closed
subsets of $2^\IN$ in computable analysis~\cite{BP03,IK21}. We denote these
spaces by $\AA$, $\AA_+$ and $\AA_-$, depending on whether we use full,
positive or negative information, respectively. 
More precise definitions follow in the next section. 

Altogether we will study three different general versions of PAC learning.
Besides the usual concepts of proper and improper PAC learning, it turned
out to be useful to consider also a variant of PAC learning that
we refer to as {\em relative PAC learning}, where besides the concept class
$\CC$ we also provide the larger class $\HH\supseteq\CC$ of hypotheses that
are permitted as results. In this case we consider the class $\CC$ given 
by positive information and the permitted hypotheses are described negatively,
as this allows one to exclude unwanted hypotheses.

As usual in Weihrauch complexity~\cite{BGP21}, we use partial multivalued
maps $f:\In X\mto Y$ to express {\em mathematical problems}. The idea is that
$\dom(f)$ is the set of {\em instances of the problem}, whereas $f(x)\In Y$
is the set of {\em solutions of the problem} for a particular instance $x\in\dom(f)$.
In this way we can express our three variants of the PAC learning problem as follows.

\begin{definition}[PAC learning problems]
	\label{def:PAC-problems}
	We define the following.
	\begin{enumerate}
		\item The {\em relative PAC learning problem}
		\begin{eqnarray*}
			\RPAC:\In\AA_+\times\AA_-\times\IN&\mto&\LL\times\IN^\IN\\ (\CC,\HH,d)&\mapsto&\{(A,m):A\mbox{ PAC learns }\CC\mbox{ relative to }\HH\mbox{ with }m\}
		\end{eqnarray*}
		with $\dom(\RPAC)=\{(\CC,\HH,d):\VCdim(\CC)\leq d\mbox{ and }\varnothing\not=\CC\In\HH\}$.
		\item The {\em proper PAC learning problem}
		\begin{eqnarray*}
			\PPAC:\In \AA \times\IN&\mto&\LL\times\IN^\IN\\
			(\CC,d)&\mapsto&\{(A,m):A\mbox{ properly PAC learns }\CC\mbox{ with }m\}
		\end{eqnarray*}
		with $\dom(\PPAC)=\{(\CC,d):\VCdim(\CC)\leq d$ and $\CC\neq\varnothing\}$.
		\item The {\em improper PAC learning problem}
		\begin{eqnarray*}
			\IPAC:\In \AA \times\IN&\mto&\LL\times\IN^\IN\\ 
			(\CC,d)&\mapsto&\{(A,m):A\mbox{ improperly PAC learns }\CC\mbox{ with }m\}
		\end{eqnarray*}
		with $\dom(\IPAC)=\{(\CC,d):\VCdim(\CC)\leq d$ and $\CC\neq\varnothing\}$.
	\end{enumerate}
	We use analogous notations $\PPAC^\pm$ and $\IPAC^\pm$ 
	with the upper index $+$ or $-$ if we replace
	the space $\AA$ by either $\AA_+$ or $\AA_-$, respectively.
\end{definition}

Regarding question (B) we will see that providing an upper bound on the VC dimension
makes a difference in some cases and no difference in other cases.
In order to be able to express these differences precisely, we introduce
variants of all the above learning problems with a fixed bound on 
the VC dimension that is no longer provided as input information.

\begin{definition}[Bounded PAC learning problems]
	We define problems 
    $\RPACd{d}$,
	$\PPACd{d}$
	$\IPACd{d}$
	for all $d\in\IN_\infty$ as in Definition~\ref{def:PAC-problems},
	but without the input $d\in\IN$ and using the domains:
\begin{enumerate}
	\item $\dom(\RPACd{d})=\{(\CC,\HH):\VCdim(\CC)<d$ and $\varnothing\not=\CC\In\HH\}$.
	\item $\dom(\PPACd{d})=\{\CC:\VCdim(\CC)<d$ and $\CC\neq\varnothing\}$.
	\item $\dom(\IPACd{d})=\{\CC:\VCdim(\CC)<d$ and $\CC\neq\varnothing\}$.
\end{enumerate}
We use analogous notations $\PPACd{d}^\pm$ and $\IPACd{d}^\pm$ 
with the upper index $+$ or $-$ if we replace
the space $\AA$ by either $\AA_+$ or $\AA_-$, respectively.
\end{definition}

For proper and improper PAC learning we can summarize some of our results in 
the table in Figure~\ref{fig:table-PAC}. The given notions of computability
are in each case optimal upper bounds with respect to the types of
notions that are mentioned in the table.

\begin{figure}[htb]
	\begin{small}
		\begin{tabular}{ccc}
			{\bf input} & {\bf proper PAC learning} & {\bf improper PAC learning}\\[0.1cm]\hline
			&&\\[-0.3cm]
			$\CC\in\AA, d>\VCdim(\CC)$ & computable & computable\\
			$\CC\in\AA_+, d>\VCdim(\CC)$ & limit computable & non-deterministically computable\\
			$\CC\in\AA_-, d>\VCdim(\CC)$ & limit computable & computable\\[0.1cm]\hline
			&&\\[-0.3cm]
			$\CC\in\AA$   & finite mind-change computable & finite mind-change computable\\
			$\CC\in\AA_+$ & limit computable & limit computable\\
			$\CC\in\AA_-$  & double limit computable & double limit computable
		\end{tabular}
	\end{small}
	\caption{Computability properties of $\PPAC,\PPACd{\infty}$ and $\IPAC,\IPACd{\infty}$, respectively.}
	\label{fig:table-PAC}
\end{figure}

However, using Weihrauch complexity we can get a much more nuanced picture
compared to the one that is expressed in the table in Figure~\ref{fig:table-PAC}.
We can also answer questions of completeness for the respective classes and
give more precise bounds (see also the diagram in Figure~\ref{fig:costs}). 
We can, for instance, also derive conclusion as the
following: improper PAC learning from positive information with given upper bound
$\IPAC^+$ is not complete for non-deterministic computations, but yet there
is no probabilistic algorithm for it and, in particular it is not Las Vegas computable in the sense of~\cite{BGH15a}.
This is because there are computable instance for which all solutions are of
PA Turing degree.
The diagram in Figure~\ref{fig:Benchmark} shows some of the non-computable
versions of PAC learning (in red) together with certain benchmark problems from
Weihrauch complexity (in blue, all problems will be defined below).

\begin{figure}[htb]
	\begin{center}
		\begin{tikzpicture}[scale=0.9,auto=left,thick,every node/.style={fill=blue!20},
				PAC/.style ={fill=violet!20}]
			\node (CNP) at (4,8) {$\C_\IN'$};
			\node (SORTP) at (4,10) {$\SORT'$};
			\node (limCNP) at (0,10) {$\lim\times\C_\IN'$};
			\node (lim) at (0,8) {$\lim$};
			\node (CR) at (0,6) {$\C_\IR\equivSW\C_\IN\times\C_{2^\IN}$};
			\node  (SORT) at (4,6) {$\SORT\equivSW\sup$};
			\node (WKL) at (-2,4) {$\C_{2^\IN}\equivSW\WKL$};
			\node (limD) at (2,4) {$\lim_\Delta\equivSW\id\times\C_\IN$};
			\node (DNCS) at (-4,2) {$\DNC_*$};
			\node (WWKL) at (0,2) {$\WWKL$};
			\node (CN) at (4,2) {$\C_{\IN}\equivSW\max$};
			\node (PA) at (-6,0) {$\PA$};
			\node (DNC) at (-2,0) {$\DNC_\IN$};
			\node (LPO) at (6,0) {$\LPO$};
			\node  (LLPO) at (2,0) {$\LLPO\equivSW\C_2$};
			\node[PAC] (PPACMI) at (-5,10) {$\PPACd{\infty}^-$};			
			\node[PAC] (IPACP) at (-3,3) {$\IPAC^+$};
			\node[PAC] (PPAC) at (-5,8) {$\PPAC^\pm\equivSW\PPACd{\infty}^+$};
			\node[PAC] (PPACI) at (8,4) {$\PPACd{\infty}$};
			\node[PAC] (RPACI) at (-5,6) {$\RPACd{\infty}$};			
			\node[PAC] (RPAC) at (-5,4) {$\RPAC$};
			\node[PAC] (VCDIM) at (8,6) {$\VCdim^+$};	
			\node[PAC] (VCDIMI) at (8,2) {$\VCdim^+_{<\infty}$};	
			\node[PAC] (VCDIMM) at (8,10) {$\VCdim^-$};	
			\node[PAC] (VCDIMMI) at (8,8) {$\VCdim^-_{<\infty}$};				
			\draw[->] (SORTP) edge (CNP);
			\draw[->] (CNP) edge (SORT);			
			\draw[->] (limCNP) edge (lim);
			\draw[->] (limCNP) edge (CNP);
			\draw[->] (lim) edge (CR);
			\draw[->] (CR) edge (WKL);
			\draw[->] (lim) edge (SORT);
			\draw[->] (SORT) edge (CN);
			\draw[->] (CR) edge (limD);
			\draw[->] (limD) edge (CN);			
			\draw[->] (LPO) edge (LLPO);
			\draw[->] (WKL) edge (WWKL);
			\draw[->] (WKL) edge (IPACP);
			\draw[->] (IPACP) edge (DNCS);
			\draw[->] (WWKL) edge (LLPO);
			\draw[->] (CN) edge (LPO);
			\draw[->] (DNCS) edge (PA);
			\draw[->] (DNCS) edge (DNC);
			\draw[->] (WWKL) edge (DNC);
			\draw[->] (WWKL) edge (LLPO);
			\draw[<->] (PPACMI) edge (limCNP);
			\draw[<->] (PPAC) edge (lim);
			\draw[<->] (RPAC) edge (WKL);
			\draw[<->] (RPACI) edge (CR);
			\draw[<->] (PPACI) edge (limD);			
			\draw[<->] (VCDIM) edge (SORT);
			\draw[<->] (VCDIMI) edge (CN);
			\draw[<->] (VCDIMM) edge (SORTP);
			\draw[<->] (VCDIMMI) edge (CNP);			
		\end{tikzpicture}
		\caption{PAC learning in the strong Weihrauch lattice.}
		\label{fig:Benchmark}
	\end{center}
\end{figure}

In the following Section~\ref{sec:Weihrauch} we provide the necessary
definitions from computable analysis and Weihrauch complexity that we
are going to use for our study. In Section~\ref{sec:VC} we study
the classification of VC dimension, which is of independent interest.
The complexity of computing the VC dimension restricted to sets
of finite VC dimension will play a role in our classification of
PAC learning problems, in particular in the cases where no upper bound
on the VC dimension is provided as input.
Our results also directly lead to a classification of the Borel complexity 
of PAC learnability.
In Section~\ref{sec:RPAC} we then study relative PAC learning
and we prove that it is equivalent to Weak K\H{o}nig's Lemma $\WKL$.
In Section~\ref{sec:PPAC} we study proper PAC learning and we prove 
that it is equivalent to the limit operation on Baire space,
no matter whether concept classes are represented positively or negatively. 
Finally, in Section \ref{sec:IPAC} we study improper PAC learning from positive information
and we show that it is strictly in between the problem $\PA$ (of computing
complete extensions of Peano arithmetic) and $\WKL$ (Weak K\H{o}nig's Lemma).
This classification allows us to draw the conclusions that we have already
mentioned above.

For the study of improper PAC learning it turned out to be
very useful to use the concept of a {\em witness function} that
was introduced by Delle Rose, Kozachinskiy, Rojas, and Steifer~\cite{DRKRS23}.

\begin{definition}[Witnesses]
	Let $\CC\In2^\IN$ be a concept class. 
	A function $f:\IN^{k+1}\to\{0,1\}^{k+1}$ is called 
	a {\em $k$--witness} or just a {\em witness function} for $\CC$, if
	\[f(x_1,...,x_{k+1})\not=(h(x_1),...,h(x_{k+1}))\]
	for all $x_1<...<x_{k+1}\in\IN$ and $h\in\CC$.
\end{definition}

Using this definition we can characterize the notion of VC dimension by
\[\VCdim(\CC)=\min\{k\in\IN:\CC\mbox{ admits a $k$--witness}\},\]
with the understanding that $\min(\varnothing)=\infty$.
The notion of a $k$--witnesses was used by the aforementioned authors to
define the notion of {\em effective VC dimension}
\[\eVCdim(\CC)=\min\{k\in\IN:\CC\mbox{ admits a computable $k$--witness}\}.\]
Using this concept, the authors proved the following non-uniformly computable version
of the Fundamental Theorem of Statistical Learning,
building on earlier results of Sterkenburg~\cite{Ste22a}, who implicitly 
proved one direction of this equivalence.

\begin{theorem}[Delle Rose, Kozachinskiy, Rojas, and Steifer 2023]
	\label{thm:Fundamental-computable}
	A concept class $\CC\In2^\IN$ admits an improper computable learner $A:\SS\to2^\IN$
	(with a computable sample complexity)
	if and only if its effective VC dimension is finite.
\end{theorem}

We start Section~\ref{sec:IPAC} with an analysis of the uniform computational 
content of the proof of 
this result, which will help us to locate the complexity of $\IPAC^+$.

Finally, in the conclusions in Section~\ref{sec:conclusion} we review our results
an we formulate a number of open questions.

\section{Preliminaries on Computable Analysis and Weihrauch Complexity}
\label{sec:Weihrauch}

We introduce some concepts from computable analysis and Weihrauch complexity
and we refer the reader to \cite{BH21,Wei00} for further details.
We follow the representation based approach to computable analysis and we recall
that a {\em representation} of a space $X$ is a surjective partial map $\delta_X:\In\IN^\IN\to X$.
In this case $(X,\delta_X)$ is called a {\em represented space}.
If we have two represented spaces $(X,\delta_X)$ and $(Y,\delta_Y)$, then
a partial function $F:\In\IN^\IN\to\IN^\IN$
is called a {\em realizer} of some partial multivalued function $f:\In X\mto Y$, if
\[\delta_YF(p)\in f\delta_X(p)\]
for all $p\in\dom(f\delta_X)$. In this situation we also write $F\vdash f$.
Any partial multivalued function $f:\In X\mto Y$ on represented spaces
is called a {\em problem} and a problem is called {\em computable}
if and only if it has a computable realizer. Other notions, such as continuity,
finite mind-change computability or limit computability can be defined for problems
analogously.

A {\em computable metric space} $(X,d,\alpha)$ is a space such that $d:X\times X\to\IR$ is a metric on $X$ with a dense sequence $\alpha:\IN\to X$ such that $d\circ(\alpha\times\alpha):\IN^2\to\IR$ is computable. Such spaces are represented by their {\em Cauchy representation} $\delta_X$,
which can be defined by
\[\delta_X(p)=x:\iff\lim_{n\to\infty}\alpha p(n)=x\mbox{ and }(\forall j)(\forall i>j)\;d(\alpha p(i),\alpha p(j))<2^{-j}.\]
The {\em real numbers} $\IR$ are represented as a computable metric space
with a standard enumeration $q:\IN\to\IQ$ of the rational numbers.
Sometimes we also use one-sided representations of the reals.
For instance, $\overline{\IR}_<=\IR\cup\{\infty\}$ is represented by
\[\delta_{\overline{\IR}_<}(p):=\sup_{n\in\IN}q_{p(n)}\]
and by $\IR_<$ we denote the restriction of this space to the reals without infinity.
The space $\IR_>$ is defined symmetrically. For every represented space $X$
there is also a canonical representation of the {\em sequence space} $X^\IN$
and of the {\em space of finite words} $X^*$ over $X$.

By $\IN_{\infty}=\IN\cup\{\infty\}$ we denote the space of {\em conatural numbers},
represented by $\delta_{\IN_\infty}(0^n1p):=n$ for $n\in\IN$ and $p\in\IN^\IN$
and $\delta_{\IN_\infty}(\widehat{0})=\infty$. This space is computably isomorphic
to the compact subspace $\{0\}\cup\{2^{-n}:n\in\IN\}$ of $\IR$.

For a computable metric space $(X,d,\alpha)$ we denote by $B(x,\varepsilon):=\{y\in X:d(x,y)<\varepsilon\}$  the {\em open ball} around $x\in X$ with radius
$\varepsilon>0$ and by $B_{\langle n,k\rangle}:=B(\alpha(n),q_k)$ we obtain a numbering of the open rational balls,
where $\langle n,k\rangle\in\IN$ denotes the {\em Cantor pairing} of the natural numbers $n,k\in\IN$.
In the case of the Cantor space $2^{\IN}$, we usually consider a bijective computable enumeration $\w:\IN\to\{0,1\}^*,i\mapsto \w_i$ of the finite words. In this case we also use the open balls $B_n:=\w_n 2^{\IN}$.

By $\OO(X)$ and $\AA(X)$ we denote the {\em spaces of open} and of {\em closed subsets} of $X$,
respectively. In fact, we need three different versions of the space of closed subsets
with respect to positive, negative and full information that we denote by
$\AA_+(X)$, $\AA_-(X)$ and $\AA(X)$, respectively.
These spaces are defined as represented spaces.
We write $\range(p)-1:=\{n\in\IN:n+1\in\range (p)\}$ for $p\in\IN^\IN$ in order to allow for the constant
zero sequence as a name of the empty set.

\begin{definition}[Hyperspace representations]
	Let $X$ be a computable metric space. We define the following representations with $p,q\in\IN^\IN$ and closed $A\In X$:
	\begin{enumerate}
		\item $\delta_{\OO(X)}$ of $\OO(X)$ by
		      $\delta_{\OO(X)}(p):=\bigcup_{n\in\range(p)-1} B_n$.
		\item $\delta_{\AA_-(X)}$ of $\AA_-(X)$ by
		      $\delta_{\AA_-(X)}(p):=X\setminus\delta_{\OO(X)}(p)$.
		\item $\delta_{\AA_+(X)}$ of $\AA_+(X)$ by
		      $\delta_{\AA_+(X)}(p)=A:\iff \range(p)-1=\{n\in\IN:B_n\cap A\not=\varnothing\}$.
		\item $\delta_{\AA(X)}$ of $\AA(X)$ by $\delta_{\AA(X)}\langle p,q\rangle=A:\iff\delta_{\AA_+(X)}(p)=\delta_{\AA_-(X)}(q)=A$.
	\end{enumerate}
\end{definition}

Here $\langle.\rangle$ denotes some
standard pairing function on Baire space $\IN^\IN$ (we use this notation for pairs as well as for the pairing of sequences).
The above representations of hyperspaces are widely used in computable analysis~\cite{BP03,Wei00,BH21}.

A point $x\in X$ in a represented space $(X,\delta)$ is called {\em computable}, if it has a computable
name, i.e., if there is a computable {\em name} $p\in\IN^\IN$ with $\delta(p)=x$. 
The representation $\delta$
also induces a topology on $X$, namely its {\em final topology}, which is the largest topology on $X$
that makes $\delta$ continuous. For the introduced hyperspaces we obtain the following (see~\cite{BP03}):
\begin{enumerate}
	\item The computable points in the space $\AA_+(X)$ are called {\em c.e.\ closed}.
	      The induced topology on the space $\AA_+(X)$ is the {\em lower Fell topology}.
	\item The computable points in the space $\AA_-(X)$ are called {\em co-c.e.\ closed}.
	      The induced topology on the space $\AA_-(X)$ is the {\em upper Fell topology}.
	\item The computable points in the space $\AA(X)$ are called {\em computably closed}.
	      The induced topology on the space $\AA(X)$ is the {\em Fell topology}.
	\item The computable points in the space $\OO(X)$ are called {\em c.e.\ open}.
	      The induced topology on the space $\OO(X)$ is the {\em Scott topology}.
\end{enumerate}
This terminology for computable sets generalizes the classical notions of c.e., co-c.e.\ and computable sets $A\In\IN$,
as these are exactly the computable points in $\AA_+(\IN)$, $\AA_-(\IN)$, and $\AA(\IN)$, respectively.
As we will mostly use these notations for Cantor space $X=2^\IN$, we also briefly write
$\AA_+$, $\AA_-$ and $\AA$ in this special case.

Occasionally, it is useful for us to use another characterization of the space $\AA_+(X)$.
Namely, if $X$ is a complete metric space, then we can also encode a closed subset $A\In X$
by enumerating a dense subset~\cite{BP03}. By $\overline{A}$ we denote the {\em closure} of a subset $A\In X$.

\begin{proposition}[Sequential representation of closed subsets]
	\label{prop:seq-positive}
	Let $X$ be a computable metric space. Then
	\[S:X^\IN\to\AA_+(X),(x_n)_{n\in\IN}\mapsto\overline{\{x_n:n\in\IN\}}\]
	is computable and if $X$, additionally, is complete, then $S$ admits
	a multivalued computable inverse, defined for all non-empty $A\in\AA_+(X)$.
\end{proposition}

It is also well-known that the closed subsets of Cantor space with full information can be identified with pruned trees~\cite{BP03}.
By $\Tr$ we denote the space of binary trees represented by their characteristic
functions.

\begin{proposition}[Pruned trees]
	\label{prop:trees}
	The map 
	\[T:\AA\to\Tr, A\mapsto T_A:=\{w\in\{0,1\}^*:w2^\IN\cap A\not=\varnothing\}\] is a computable embedding, i.e., it is computable, injective and has a computable partial inverse.
\end{proposition}

We note that if $[T]:=\{p\in 2^\IN:(\forall n\in\IN)\;p|_n\in T\}$ denotes
the {\em set of infinite paths} of a tree $T\in\Tr$, then $[T_A]=A$ for all $A\in\AA$.

At this point we briefly discuss how the effectivity notions that correspond
to considering hypothesis classes $\HH\in\AA_+$ relate to 
other notions that have been used in the literature. Agarwal et al.~\cite{AAB+20}
introduced
so-called {\em RER hypothesis classes} 
$\HH_e=\{\varphi_i:i\in\W_e\}$,
where $\varphi$ is a standard G\"odel numbering of the partial computable
functions $\varphi_i:\In\IN\to\IN$ and $\W_e:=\dom(\varphi_e)$ is the derived
standard numbering of the computably enumerable sets $W_e\In\IN$.
Here $\HH_e$ is only defined for such $e\in\IN$ for which all the functions $\varphi_i$
with $i\in\W_e$ are total. A related concept of {\em DR hypothesis classes} is defined
via decidable sets of indices instead of the enumerations. 
By Proposition~\ref{prop:seq-positive} our c.e.\ closed hypothesis classes
$\HH$ are exactly the closures of RER classes. However, when it comes to
uniform computability, then using G\"odel numberings is less natural, because
the function $f:\In\IN\to\AA_+,e\mapsto\overline{\HH_e}$ is computable,
but it has no computable left inverse, not even restricted to singletons $\HH_e$.
This is because finding a G\"odel number for a given computable sequence is a non-trivial
computational problem by itself (whose Weihrauch complexity has been studied in~\cite{Bra23a}). 
In other words, working with G\"odel numbers artificially blurs the complexity of PAC learning
when one looks at it from a uniform perspective and this is the reason why
we have chosen an approach that is better suited for this purpose.
It has already been noticed by Akbari and Harrison-Trainor~\cite{AH25} that
working with G\"odel numbers can be unnatural in the context
of PAC learning and the solution proposed by them is to consider 
{\em relativized RER classes}
$\HH_{\langle p,e\rangle}=\{\varphi_i^p:i\in\W_e^p\}$,
where everything is relativized to some oracle $p\in\IN^\IN$.
Again $\HH_{\langle p,e\rangle}$ is defined just for such $\langle p,e\rangle$
for which the class contains only total functions.
While in \cite{AH25} these notions are considered on a cone, we basically
work directly with closures of such classes. We make this statement, which
is a consequence of Proposition~\ref{prop:seq-positive}, more precise.

\begin{corollary}[Relativized RER classes]
	The function $F:\In\IN^\IN\to\AA_+,\langle p,e\rangle\to\overline{\HH_{\langle p,e\rangle}}$ is computable, surjective and admits a computable left inverse.
\end{corollary}

This corollary might help the reader to translate our results for the space $\AA_+$
into the setting of RER classes, if desired. In many of our constructions 
the classes are already closed relativized RER classes (in those cases no 
closures need to be taken). 

We now proceed to introduce concepts from Weihrauch complexity (see \cite{BGP21} for a survey). Here $\id:\IN^\IN\to\IN^\IN$ denotes the identity on Baire space.

\begin{definition}[Weihrauch reducibility]
	Let $f:\In X\mto Y$ and $g:\In Z\mto W$ be problems. We say that
	\begin{enumerate}
		\item $f$ is {\em Weihrauch reducible} to $g$, in symbols $f\leqW g$, if there are computable functions
		      $H,K:\In\IN^\IN\to\IN^\IN$ such that $H\langle\id,GK\rangle\vdash f$, whenever $G\vdash g$ holds.
		\item $f$ is {\em strongly Weihrauch reducible} to $g$, in symbols $f\leqSW g$, if there are computable
		      $H,K:\In\IN^\IN\to\IN^\IN$ such that $HGK\vdash f$, whenever $G\vdash g$ holds.
	\end{enumerate}
\end{definition}

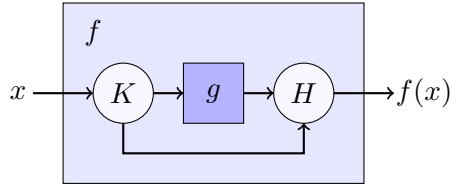
\begin{figure}[htb]
	\begin{tikzpicture}[scale=.4,auto=left]
		\draw[style={fill=blue!10}] (-3,6) rectangle (7,0);
		\draw[style={fill=blue!30}]  (1,4) rectangle (3,2);
		\draw[style={fill=blue!03}]  (-1,3) ellipse (1 and 1);
		\draw[style={fill=blue!03}]  (5,3) ellipse (1 and 1);
		\node at (-1,3) {$K$};
		\node at (5,3) {$H$};
		\node at (2,3) {$g$};
		\node at (-2,5) {$f$};
		\node at (-4.5,3) {$x$};
		\node at (9,3) {$f(x)$};
		\draw[->,thick] (0,3) -- (1,3);
		\draw[->,thick] (-4,3) -- (-2,3);
		\draw[->,thick] (3,3) -- (4,3) ;
		\draw[->,thick] (6,3) -- (8,3);
		\draw[->,thick] (-1,2) -- (-1,1) -- (5,1) -- (5,2);
	\end{tikzpicture}
\caption{Weihrauch reducibility $f\leqW g$.}
\label{fig:Weihrauch}
\end{figure}

We use the notations $\equivW$ and $\equivSW$ for the equivalence relations induced by the
preorders $\leqW$ and $\leqSW$, respectively. The corresponding equivalence classes
are called {\em Weihrauch degrees} and {\em strong Weihrauch degrees}, respectively.
There is also a purely topological version of (strong) Weihrauch reducibility denoted by
$\leq_{\mathrm W}^*$ and $\leq_{\mathrm{sW}}^*$, respectively. In this case $H,K$
are just required to be continuous. All our results also hold true for continuous
Weihrauch reducibility. The diagram in Figure~\ref{fig:Weihrauch} illustrates the reducibility.

An important class of problems that can be used to calibrate other problems are the
choice problems. The {\em choice problem} $\C_X$ of a computable metric space $X$
is the problems of finding a point in a negatively represented closed set $A\in\AA_-(X)$.

\begin{definition}[Choice problem]
	The {\em choice problem} of a computable metric space $X$ is defined by
	\[\C_X:\In\AA_-(X)\mto X,A\mapsto A\]
	with $\dom(\C_X):=\{A\in\AA_-(X):A\not=\varnothing\}$.
\end{definition}
We will use the choice problems $\C_2$, $\C_\IN$, $\C_{2^\IN}$ and $\C_\IR$.
Choice for the two point space $2=\{0,1\}$ is also known as $\LLPO$, the
	{\em lesser limited principle of omniscience}.
The {\em limited principle of omniscience} $\LPO:\IN^\IN\to\{0,1\}$ is the
characteristic function of the singleton set $\{\widehat{0}\}$ with the constant zero sequence
$\widehat{0}$. In general, $\widehat{n}\in\IN^\IN$ denotes the constant sequence with value $n\in\IN$.

We define a number of further benchmark problems that we are going to use.
We often identify the numbers $n\in\IN$ with $X=\{0,...,n-1\}$.
The limit in the following definition is understood with respect to the product
topology on $\IN^\IN$. We use the relation $\ll$ of relative PA degrees,
i.e., $p\ll q$ for $p,q\in\IN^\IN$ denotes the property that every
infinite binary $p$--computable tree has a $q$--computable path.
In this case one says that $q$ is has {\em PA degree relative to} $p$.
The terminology originates from the fact that for computable $p$ this
corresponds to the complexity of finding a complete extension of Peano arithmetic.

\begin{definition}[Problems]
We consider the following problems for all $X\in\IN\cup\{\IN\}$:
\begin{enumerate}
\item $\lim:\In\IN^\IN\to\IN^\IN,\langle p_0,p_1,p_2,...\rangle\mapsto\lim_{i\to\infty}p_i$
is called the {\em limit problem}. 
\item $\lim_\Delta:\In\IN^\IN\to\IN^\IN$ is the restriction of $\lim$ to eventually constant sequences $(p_i)$.
\item $\WKL:\In\Tr\mto2^\IN,T\mapsto[T]$ is called {\em Weak K\H{o}nig's Lemma}.
\item $\WWKL$ is the restriction of $\WKL$ to trees $T$ with $[T]$ 
of positive measure and is called {\em Weak Weak K\H{o}nig's Lemma}. 
\item $\SORT: 2^{\IN}\to 2^{\IN},p\mapsto\widehat{0}$ if $p\in 2^{\IN}$ contains infinitely many zeroes, and $p\mapsto 0^n \widehat{1}$ if $p$ contains exactly $n$ zeroes
is called the {\em sorting problem}.
\item $\sup: \IN^\IN \to \IN_\infty,p\mapsto\sup \{ p(i):i \in\IN \}$ is called
the {\em supremum problem}.
\item $\max:\In\IN^\IN\to\IN,p\mapsto\max\{p(i):i\in\IN\}$ 
is called the the {\em maximum problem}.
\item $\B:\In\IN^\IN\mto\IN,p\mapsto\{n\in\IN:(\forall k\in\IN)\;p(k)\leq n\}$
is called the {\em boundedness problem}.
\item $\ACC_X:\In\AA_-(X)\mto X$ is the restriction of $\C_X$ to sets
$A\in\AA_-(X)$ with $|X\setminus A|\leq 1$ and is called {\em all or co-unique choice problem}.
\item $\DNC_X:\IN^\IN\to X^\IN,p\mapsto\{q:(\forall n)\;\varphi_n^p(n)\not=q(n)\}$ is called
the {\em diagonally non-computability problem}.
\item $\PA:\IN^\IN\to\IN^\IN,p\mapsto\{q:p\ll q\}$ is the problem of {\em degrees of Peano arithmetic}.
\end{enumerate}
All problems are understood with their natural domains, i.e., $\lim$ is defined for 
converging sequences, $\WKL$ for infinite trees, $max$ and $\B$ for bounded sequences, etc.
\end{definition}

It basically follows from results in \cite{BHK17a,BHK17} that the given problems
are ordered as shown in Figure~\ref{fig:Benchmark} (where an arrow $f\leftarrow g$ means $f\leqSW g$).

We also need a number of operators on problems that are commonly used in Weihrauch complexity.
The {\em jump} of a represented space $(X,\delta)$ is the represented space $(X,\delta\circ\lim)$ and
we denote it also by $X'$.

\begin{definition}[Operations on problems]
	Let  $f:\In X\mto Y$, $g:\In Z\mto W$ be problems. We define the following operations:
	\begin{enumerate}
		\item The {\em parallelization}
		      $\widehat{f}:\In X^\IN\mto Y^\IN,(x_n)_{n\in\IN}\mapsto\bigtimes_{n\in\IN} f(x_n)$.
		\item The {\em product}
		      $f\times g:\In X\times Z\mto Y\times W,(x,z)\mapsto f(x)\times g(z)$.
		\item The {\em compositional product}
		      $f*g:=\max\nolimits_{\leqW}\{F\circ G:F\leqW f\mbox{ and }G\leqW g\}$.
            \item The {\em jump} $f':\In X'\mto Y,x\mapsto f(x)$.
	\end{enumerate}
\end{definition}

We note that all these operations are monotone with respect to Weihrauch reducibility and hence
they extend to Weihrauch degrees, except for the jump, which is monotone only with respect to 
strong Weihrauch reducibility~\cite{BGM12}.
Intuitively, the compositional product captures the most powerful problem that
can be obtained by first using $g$ and then (after some intermediate computation) $f$.
The compositional product was introduced in~\cite{BGM12}. That the corresponding
maximum always
exists was proved in~\cite{BP18}. Strictly speaking, $f*g$ only defines an ordinary Weihrauch
degree, but there is a specific instance of this degree that is maximal with respect
to strong Weihrauch reducibility and that also allows us
to formulate strong Weihrauch reductions involving compositional products (see~\cite{BP18,BGP21} for details).

We will need a number of characterizations of the degree of
Weak K\H{o}nig's Lemma $\WKL$ that were essentially proved in~\cite{BG11,BBP12,BHK17a}.
We recall that $\LL=(2^\IN)^\SS$ denotes the set of learners that is computably
isomorphic to Cantor space $2^\IN$ via some standard coding of $\IN\times\SS$.

\begin{proposition}[Weak K\H{o}nig's Lemma]
	\label{prop:WKL}
	$\WKL\equivSW\C_{2^\IN}\equivSW\C_\LL\equivSW\widehat{\C_2}\equivSW\DNC_2$.
\end{proposition}

Choice on the natural numbers $\C_\IN$ can be characterized as follows~\cite{BG11a,BBP12,BGP21}.

\begin{proposition}[Choice on the natural numbers]
	\label{prop:CN}
	$\B\equivW\C_\IN\equivSW\max\equivW\lim_\Delta\equivSW\id\times\C_\IN$
and $\B\lSW\C_\IN\lSW\lim_\Delta$.
\end{proposition}

Choice on Euclidean space $\C_\IR$ can be written in different ways as
product of $\C_\IN$ and $\C_{2^\IN}$, namely $\C_\IR\equivSW\C_\IN\times\C_{2^\IN}\equivSW\C_{2^\IN}*\C_\IN$,
which was essentially proved in~\cite{BBP12,BGM12}.
Even though $\B\lSW\C_\IN$, we obtain $\B\times\C_{2^\IN}\equivSW\C_\IN\times\C_{2^\IN}$,
because the input can be copied through $\C_{2^\IN}$ (as $\id\times\C_{2^\IN}\leqSW\C_{2^\IN}$).

\begin{proposition}[Choice on Euclidean space]
	\label{prop:C_R}
	$\C_\IR\equivSW\C_{2^\IN}*\C_\IN\equivSW\C_\IN\times\C_{2^\IN}\equivSW\B\times\C_{2^\IN}$.
\end{proposition}

The following characterization is well-known and was proved in~\cite{BG11a,Bra18}.

\begin{proposition}[Limit]
	\label{prop:lim}
	$\lim\equivSW\widehat{\LPO}\equivSW\lim*\C_\IN$.
\end{proposition}

The following characterization was first proved by Higuchi and Kihara~\cite[Proposition~7.9]{HK14} (see also \cite[Theorem~5.2]{BHK17a}).
We note that $\ACC_n$ is also known under the name $\LLPO_n$.

\begin{proposition}[Diagonally non-computability]
\label{prop:DNC}
$\DNC_X\equivSW\widehat{\ACC_X}$ for all $X\in\IN\cup\{\IN\}$ with at least two elements.
\end{proposition}

We introduce a new variant of $\DNC_n$ that is relevant for the analysis of 
improper PAC learning and that we informally define by $\DNC_*:=\sqcap_{k\geq2}\DNC_k$. 
The idea is that it corresponds to the logical problem $(\exists k\geq2)\;\DNC_k$.
As we did not define the operation $\sqcap$ for sequences of problems in general\footnote{For a sequence of problems $(f_k)_{k\in\IN}$ we can define $\sqcap_{k\in\IN}f_k$ (see \cite{BGH15a})
such that it is the infimum of the $f_i$ with respect to continuous Weihrauch
reducibility. However, with respect to computable Weihrauch reducibility
there are no non-trivial infima of natural problems~\cite[Corollary~3.18]{HP13}.},
we define the problem $\DNC_*$ ad hoc here.

\begin{definition}[Finitary DNC]
We define $\DNC_*:\In\bigtimes_{k=2}^\infty\AA_-(k)^\IN\mto\IN^\IN$ by
\[\DNC_*(((A_{n,k})_{n\in\IN})_{k\geq2}):=\{\langle k,p\rangle:(\forall n)\;p(n)\in A_{n,k}\}\]
where $|k\setminus A_{n,k}|\leq 1$ for all $n,k$. We call this the {\em finitary DNC problem}.
\end{definition}

In other words, the instances of $\DNC_*$ are sequences $(B_k)_{k\geq 2}$,
where every $B_k=(A_{n,k})_{n\in\IN}$ is an instance of $\widehat{\ACC_k}$ and the solution 
is a pair $(k,p)$ such that $p\in\widehat{\ACC_k}(B_k)$.
We will use the problem $\DNC_*$ in order to separate improper PAC learning from
$\DNC_\IN,\PA$ and $\WWKL$. This is possible because we have the following result.

\begin{proposition}[Finitary DNC]
	\label{prop:finitary-DNC}
	$\PA\lW\DNC_*$, $\DNC_\IN\lW\DNC_*$ and $\DNC_*\nleqW\WWKL^{(n)}$ for all $n\in\IN$.
\end{proposition}
\begin{proof}
By Proposition~\ref{prop:DNC} we obtain $\DNC_\IN\equivSW\widehat{\ACC_\IN}\leqSW\DNC_*$.
The latter holds obviously, as $f:\AA_-(\IN)\times\IN\to\AA_-(k),(A,k)\mapsto A\cap k$
is computable uniformly in $k$. The reduction $\PA\leqSW\DNC_k=\widehat{\ACC_k}$ 
is well known~\cite{BHK17a} and easy to
see, as $p\ll q$ holds if and only if every $q$ computes a function $f:\IN\to k$
that is diagonally non-computable relative to $p$. This holds uniformly in $k$
and hence we obtain $\PA\leqSW\DNC_*$. Since $\PA\nleqW\WWKL^{(n)}$ for all $n\in\IN$
by \cite[Corollary~7.11]{BHK17a} we obtain $\DNC_*\nleqW\WWKL^{(n)}$ for all $n\in\IN$.
Since $\DNC_\IN\leqSW\WWKL$ by \cite[Corollary~7.7]{BHK17a}, it follows
that $\DNC_*\nleqW\DNC_\IN$. On the other hand, also $\DNC_\IN\nleqW\PA$ is
known (not even $\ACC_\IN\leqW\PA$ holds, as $\PA$ is densely realized
in the sense of \cite{BHK17a}) and hence $\DNC_*\nleqW\PA$ follows. 
\end{proof}

This result shows that the problem $\DNC_*$ is related to other well-known problems
as displayed in Figure~\ref{fig:Benchmark}.

It is well-known that the limit map is also equivalent to converting the negative information of a closed subset of Cantor space into positive information, and vice versa. We will need the following direction of this result~\cite[Propositions~4.2, 4.5]{BG09}.

\begin{proposition}[Converting information of a closed subset]
	\label{prop:convert}
	The two identity functions $\id:\AA_+\to \AA_-,A\mapsto A$ and $\id:\AA_-\to \AA_+,A\mapsto A$ are both strongly Weihrauch reducible to $\lim$.
\end{proposition}

In Corollary~\ref{cor:conversion} we are going to strengthen this result.
We close this section by mentioning that $\C_\IN$, $\C_{2^\IN}\equivSW\WKL$, $\WWKL$ and $\lim$
characterize some natural classes of problems that all have been studied
independently~\cite{BBP12,BGH15a}.

\begin{theorem}[Classes of computable problems]
	\label{thm:computable-classes}
	Let $f$ be a problem. Then:
	\begin{enumerate}
		\item $f\leqW\C_\IN\iff f$ is computable with finitely many mind changes.
		\item $f\leqW\WKL\iff f$ is non-deterministically computable.
		\item $f\leqW\WWKL\iff f$ is Las Vegas computable.
		\item $f\leqW\lim\iff f$ is limit computable.
	\end{enumerate}
\end{theorem}

We do not define the involved notions of computability separately here, but
the characterization in Theorem~\ref{thm:computable-classes} can also be taken as a definition. We just note that
finite mind change computations and limit computations have often been applied
in the context of algorithmic learning theory~\cite{Gol65}.

\section{VC Dimension}
\label{sec:VC}

The concept of {\em VC dimension} was introduced by Vapnik and Chervonenkis~\cite{VC71}.
As it has found many applications, also outside of the theory of machine learning,
it is of independent interest to provide a precise classification of its complexity.
For our study of PAC learning we will mostly need the results that are restricted
to sets of finite VC dimension.

We will show that the complexity of computing the VC dimension of closed subsets $\CC\in \AA$ with full information or for closed subsets $\CC \in \AA_+$ with positive information is that of the binary sorting problem. For closed subsets $\CC \in \AA_-$ with negative information the problem is equivalent to the jump of the binary sorting problem.

\begin{definition}[VC dimension problem]\ 
	\begin{enumerate}
		\item By $\VCdim: \AA\to \IN_\infty$ we denote the problem of determining the VC dimension
		of closed sets, just as defined in the introduction.
		\item By $\VCdim_{<\infty}:\In\AA\to\IN$ we denote the problem of determining the VC dimension restricted to closed sets of finite VC dimension.
		\item By $\B\VCdim:\In\AA\mto\IN,\CC\mapsto\{d\in\IN:\VCdim(\CC)\leq d\}$ we denote the problem of determining an upper bound on the VC dimension of a closed set of finite VC dimension.
	\end{enumerate}
	We add an upper index $+$ or $-$ to each of these problems if it is defined on
	$\AA_+$ or $\AA_-$ instead of $\AA$, respectively.
\end{definition}

We will now use the binary sorting problem $\SORT$ as a benchmark problem.
This problem was originally introduced and studied by Neumann and Pauly~\cite{NP18}, 
even though the same ordinary Weihrauch equivalence class was already studied
by Weihrauch~\cite{Wei92a} under the name $\EC_1$ (see \cite{Bra22a} for the equivalence).
By $\SORT_{<\infty}$ we denote the sorting problem restricted to binary
sequences with only finitely many zeros.
As a preparation we prove that the binary sorting problem and the
supremum problem are actually equivalent.

\begin{lemma}
	\label{lemma:sort-benchmark}
	$\SORT\equivSW\sup$ and $\SORT_{<\infty}\equivSW\max\equivSW\C_\IN$.
\end{lemma}

\begin{proof} We prove the first equivalence.
	
	\par{$\SORT\leqSW\sup$.} Given a sequence $p\in 2^{\IN}$, produce a sequence $q\in\IN^{\IN}$ by setting $q(0)=1-p(0)$ and, for $i\in\IN$, $q(i+1)=q(i)$ if $p(i+1)=1$ and $q(i+1)=q(i)+1$ if $p(i+1)=0$. This way, $q(i)$ is the number of zeroes in $p\vert_{i+1}$. Then, apply $\sup$ to $q$ to get a $\delta_{\IN_{\infty}}$--name for the number of zeroes in $p$, which is exactly $\SORT(p)$.
	
	\par{$\sup\leqSW\SORT$.} Given a sequence $p\in\IN^{\IN}$, first write $p(0)$ zeroes and set $\texttt{sup}:=p(0)$. Then, for every read symbol $p(i)$ ($i\geq 1$), write a one, and if $p(i)> \texttt{sup}$, write $p(i)-\texttt{sup}$ zeroes and update $\texttt{sup}:=p(i)$ afterward. The resulting sequence $q$ contains infinitely many zeroes if $p$ is unbounded, and $\sup_i p(i)$ zeroes otherwise. Therefore, by applying $\SORT$ to $q$, we produce a name of $\sup(p)$.
	
	The equivalence $\SORT_{<\infty}\equivSW\max$ follows from the same construction
	above, restricted to sequences with only finitely many zeros and bounded sequences,
	respectively. One just needs to use the additional fact that $\id:\IN\into\IN_\infty$
	is a computable embedding. The equivalence $\max\equivSW\C_\IN$ is well-known~\cite[Theorem~11.7.13]{BGP21}.
\end{proof}

We can now state our main result regarding the VC dimension.
By $w\cdot p$ we denote the {\em concatenation} of a word $w\in\IN^*$ with
a sequence $p\in\IN^\IN$ (or another word $p\in\IN^*$).
By $h|_k$ we denote the {\em restriction} of a function $h:\In\IN\to\IN$ to
$k=\{0,...,k-1\}$.

\begin{theorem}[VC dimension]
	\label{thm:VC}
	$\VCdim \equivSW \VCdim^+ \equivSW \SORT$.
\end{theorem}

\begin{proof}
	By Lemma~\ref{lemma:sort-benchmark} it is sufficient to prove
	\[\SORT\leqSW\VCdim\leqSW\VCdim^+\leqSW\sup.\]
	The second reduction holds obviously. We prove the other two reductions.
	
	\begin{enumerate}
		\item \par{$\VCdim^+\leqSW \sup$.} 
		Given $\CC\in\AA_+$, we can assume by Proposition~\ref{prop:seq-positive}
		that we have a sequence $(h_n)_{n\in\IN}$ in $2^\IN$ such that
		$\CC=\overline{\{h_n:n\in\IN\}}$. We define $\CC_k:=\{h_n|_k\cdot\widehat{0}:n\leq k\}$. Then we obtain
		\[\VCdim(\CC)=\sup_{k\in\IN}\VCdim(\CC_k).\]
		This holds as any set $A\In\IN$ that is shattered by $\CC$
		must be shattered by some prefix of some $h_n$ of length $\max(A)+1$.
		The VC dimensions of $\CC_k$ can be easily computed, as only $k+1$
		finite prefixes of length $k$ need to be considered. 
		\item \par{$\SORT\leqSW \VCdim$.} Given a string $p\in2^{\IN}$, define the following subset
		\[
		\CC=\{h:\IN\to\{0,1\}: (\forall i\in\IN)(p(i)=1\implies h(i)=0)\} \in \AA.
		\]
		Observe that the VC dimension of $\CC$ is the number of indices over which the functions $h$ in $\CC$ can take both values $0$ and $1$, which by definition is exactly the numbers of zeroes in $p$.
		
		To compute $\CC\in\AA$ from $p$, we iterate over all natural numbers $i\in\IN$, starting from $i=0$. For each $i$, the extendable prefixes of length $i+1$ are those of the form $a_{0}\dots a_{i}\in \{0,1\}^{i+1}$ respecting the constraint $a_{k}=0$ for all $k\leq i$ such that $p(k)=1$. Read $p\vert_{i+1}$, identify these prefixes and write the associated balls into the positive information, then write all other prefixes of length $i+1$ into the negative information.
		
		Then, apply $\VCdim$ to get the VC dimension of $\CC$ in $\IN_\infty$. 
		This is a $\delta_{\IN_\infty}$--name of the number of zeroes in $p$, 
		which is exactly $\SORT(p)$.\qedhere
	\end{enumerate}
\end{proof}

Essentially the same proof also shows 
\[\SORT_{<\infty}\leqSW\VCdim_{<\infty}\leqSW\VCdim^+_{<\infty}\leqSW\max.\]
Hence, with the help of Lemma~\ref{lemma:sort-benchmark} we obtain the following.

\begin{corollary}
	\label{cor:VC}
	$\VCdim_{<\infty} \equivSW \VCdim^+_{<\infty} \equivSW \C_\IN$
	and $\B\VCdim\equivSW\B\VCdim^+\equivSW\B$.
\end{corollary}

The second equivalence chain follows from $\VCdim_{<\infty} \equivSW \VCdim^+_{<\infty} \equivSW \max$ using just the same reduction functions, as $\B$ just computes an upper
bound of the maximum.

We can also conclude from the proof of Theorem~\ref{thm:VC} that
the VC dimension map is lower semi-computable as a real-valued function.
This is expressed more formally in the following corollary.

\begin{corollary}
	VC dimension, considered as a map of type $\VCdim:\AA_+\to\overline{\IR}_<$, is
	computable. This also holds with $\AA$ instead of $\AA_+$.
\end{corollary}

It is interesting to note that adding the negative information does not reduce the Weihrauch complexity of computing the VC dimension.
Intuitively, the irrelevance of the negative information can be understood in the following
way: this information can help us to know that a certain subset is not shattered by the hypothesis class that we consider up to a certain depth. 
But we cannot conclude anything about the VC dimension based on that,
since a subset of the same size could very well be shattered at a greater depth of the class.

We now show what happens if we only include the negative information. Consistently with the previous result, this information cannot be exploited in a better way than translating it into positive information before using it.

\begin{theorem}[VC dimension with negative information]
	\label{thm:VCdim-negative}
	$\VCdim^-\equivSW\VCdim'\equivSW\SORT'$.
\end{theorem}

\begin{proof}
	
	$\VCdim^-\leqSW\VCdim'\leqSW \SORT'$. This follows from Theorem~\ref{thm:VC} with
	Proposition~\ref{prop:convert}.

	$\SORT'\leqSW\VCdim^-$.
	Let $(p_{n})\in(\IN^{\IN})^{\IN}$ be a sequence converging to $p\in 2^{\IN}$, and let $\ell\in\IN_\infty$ denote the number of zeroes in $p$. We want to compute the negative information of a closed set $\CC$ whose VC dimension is $\ell$. We perform the construction so that for any $h\in\CC$ and $i\in\IN$ we have $h(i)=0$, except for certain $i$'s, which we refer to as \textit{open paths}, where $h(i)$ can take both values $0$ and $1$. We want the number of open paths to be $\ell$, so that we obtain $\VCdim(\CC)=\ell$.
	
	To do that, we keep track of a variable $\mathtt{open\_paths}$ containing a list of elements of $\IN\times\IN$. For $(i,k)\in\IN\times\IN$ in this list, $i$ denotes an index which we keep open, while $k$ denotes the position of the zero in $(p_{n})$ that we associate with $i$. We also keep a variable $\mathtt{N}$, initialized by $0$, which is the latest index that we have handled in the description of $\CC$.
	
	Now, to construct $q\in\IN^{\IN}$ such that 
	\begin{eqnarray}
		\CC=2^\IN\setminus\bigcup_{m+1\in\range(q)} \w_m2^\IN=\delta_{\AA_-}(q),
	\end{eqnarray}
	where $\w:\IN\to\{0,1\}^*$ is our bijective computable numbering.
	We repeat the following steps:
	\begin{enumerate}
		\item Read $p_{\mathtt{N}\vert_{\mathtt{N+1}}}$.
		\item For any $k\leq \mathtt{N}$ so that $p_{\mathtt{N}}(k)=0$ and there is no $i\in\IN$ with $(i,k)$ in $\mathtt{open\_{paths}}$, add an open path for $k$, that is add $(\mathtt{N},k)$ to $\mathtt{open\_{paths}}$, and then increment $\mathtt{N}$.
		\item For any $k\leq\mathtt{N}$ so that $p_{\mathtt{N}}(k)\neq 0$ with some $i\in\IN$ with $(i,k)$ in $\mathtt{open\_{paths}}$ (i.e., $k$ was ``wrongly left open''), remove $(i,k)$ from $\mathtt{open\_{paths}}$.
		\item Increment $\mathtt{N}$ by one so that one additional path is closed.
		\item Effectively close all paths that have not been left open at this point, that is write all $m\in\IN$ into $q$ so that $|\w_{m-1}|\leq \mathtt{N}+1$ \textbf{and} there is $i\leq\mathtt{N}$ with $\w_{m-1}(i)=1$ and yet no $k$ such that $(i,k)$ is in $\mathtt{open\_{paths}}$.
	\end{enumerate}
	For any $k\in\IN$, there is either zero or one $i\in\IN$ such that $(i,k)$ indefinitely remains in $\mathtt{open\_{paths}}$, and there is one if and only if $p(k)=0$. Indeed, if $p(k)=0$ holds, there is a minimal $N\in\IN$ such that $p_{n}(k)=0$ for all $n\geq N$, and $k$ is added when $p_{N}(k)$ is read and then never removed. If $p(k)=1$, any addition of $k$ to $\mathtt{open\_paths}$ is eventually removed.
	
	The construction also ensures with the help of $\mathtt{N}$ 
	that no $i$ remains with different $k$'s simultaneously as pair
	$(i,k)$ in $\mathtt{open\_{paths}}$.
	Hence, there are exactly $\ell\in\IN_\infty$ values $i\in\IN$ such that 
	$(i,k)$ indefinitely remains
	in $\mathtt{open\_{paths}}$ for some $k\in\IN$. Let $i_1<i_2<...$ be these values. With $j_1:=i_1$ and $j_{k+1}:=i_{k+1}-i_k-1$ we obtain that 
	$\CC$ has the form 
	\begin{itemize}
		\item $\CC=0^{j_1}\cdot\{0,1\}\cdot0^{j_2}\cdot\{0,1\}\cdot...\cdot0^{j_\ell}\cdot\{0,1\}\cdot\widehat{0}$ if $\ell$ is finite or 
		\item $\CC=0^{j_1}\cdot\{0,1\}\cdot0^{j_2}\cdot\{0,1\}\cdot0^{j_3}\cdot\{0,1\}\cdot...$,
		if $\ell=\infty$.
	\end{itemize}
	In particular, $\VCdim(\CC)=\ell$.
	
	Therefore, once $q$ is produced, it can be fed into a realizer of $\VCdim^-$ to produce a name for the outcome of $\SORT'$.
\end{proof}

Restricted to sets of finite VC dimension and sequences with only finitely many zeros
the same consideration shows the following with the help of Corollary~\ref{cor:VC}
and Lemma~\ref{lemma:sort-benchmark}.

\begin{corollary}
	\label{cor:VCdim-negative}
	$\VCdim^-_{<\infty}\equivSW\C_\IN'$ and $\B\VCdim^-\equivSW\B'$.
\end{corollary}

We can also easily derive results on the Borel complexity of the class of PAC learnable
concept classes from our results. This result corresponds to classifications
in the arithmetical hierarchy obtained by Schaefer~\cite[Theorem~4.1]{Sch99b} (see 
also Calvert~\cite{Cal15}).
As this is just a side result for us, we do not formally discuss and introduce
the Borel hierarchy here. For the interpretation of Borel complexity for
represented spaces we refer the reader to \cite{Bra05,BG09}.

\begin{theorem}[Borel complexity of PAC learnability]
	\label{thm:Borel}
The class $\PP\LL$ of PAC learnable closed concept classes is
$\Sigma^0_2$--complete in $\AA$ and $\AA_+$, and $\Sigma^0_3$--complete in $\AA_-$
with respect to the effective Borel hierarchy.
\end{theorem}
\begin{proof}
By the Fundamental Theorem of Statistical Learning, Theorem~\ref{thm:fundamental} we have
\[\PP\LL:=\{\CC\in\AA:\CC\mbox{ PAC learnable}\}=\{\CC\in\AA:\VCdim(\CC)<\infty\mbox{ and }\CC\not=\varnothing\}.\]
The function $f:\AA\to2^\IN$ that maps every concept class
$\CC\in\AA$ to $f(\CC)=0^n\widehat{1}$ if $n=\VCdim(\CC)<\infty$ (or $\CC=\varnothing$) and
$f(\CC)=\widehat{0}$ if $n=\infty$ is strongly Weihrauch equivalent to $\SORT$ 
by Theorem~\ref{thm:VC} and hence reducible to $\lim$. This implies
that $f$ is $\Sigma^0_2$--measurable~\cite{Bra05}.
Hence
$f^{-1}(2^\IN\setminus\{\widehat{0}\})$ is an effective $\Sigma^0_2$--set
and so is $\PP\LL=f^{-1}(2^\IN\setminus\{\widehat{0}\})\setminus\{\varnothing\}$
as $\{\varnothing\}$ is co-c.e.\ closed in $\AA$.
Considering $f$ as a function of type $f:\AA_+\to2^\IN$ shows
analogously that $\PP\LL$ is also an effective $\Sigma^0_2$--set in $\AA_+$,
as $\{\varnothing\}$ is also co-c.e.\ closed in $\AA_+$.
If we consider $f$ as a function of type $f:\AA_-\to2^\IN$, then we obtain
with Theorem~\ref{thm:VCdim-negative} that $f$ is strongly Weihrauch equivalent
to $\SORT'$ and hence reducible to $\lim'$, which implies that $\PP\LL$ is
effectively $\Sigma^0_3$--measurable.
In this case we use that $\{\varnothing\}$ is c.e.\ open in $\AA_-$, due
to the compactness of Cantor space $2^\IN$.

It remains to prove the completeness claims.
It is well-known that the class $\FF\In 2^\IN$ that contains all binary sequences with only finitely many digits $1$, 
is $\Sigma^0_2$--complete in the (effective) Borel hierarchy~\cite[Exercise~23.1]{Kec95}.
We consider the computable function $g:2^\IN\to\AA$
with 
\[g:2^\IN\to\AA, p\mapsto\{0,1\}^{p(0)}\cdot\{0,1\}^{p(1)}\cdot\{0,1\}^{p(2)}\cdot...,\]
where we assume that $\{0,1\}^0=\{0\}$. Then $g(p)\not=\varnothing$ and 
\[p\in\FF\iff \VCdim(g(p))<\infty\iff g(p)\in\PP\LL\]
for all $p\in2^\IN$. 
Hence, $g$ yields an effective Wadge reduction that witnesses $\Sigma^0_2$--completeness
of $\PP\LL$ in $\AA$. This implies $\Sigma^0_2$--completeness in $\AA_+$ too.

We still need to prove $\Sigma^0_3$--completeness of $\PP\LL$ in $\AA_-$.
It is well-known that the set $\GG:=\{p\in2^\IN:(\forall^\infty k)(\exists n)\;p\langle n,k\rangle=1\}$ is (effectively) $\Sigma^0_3$--complete~\cite[Exercise~23.2]{Kec95}.
We show that there exists a computable function $h:2^\IN\to\AA_-$ such that
$p\in\GG\iff h(p)\in\PP\LL$ and such that all $h(p)$ are non-empty.
Given $p\in 2^\IN$ we compute a sequence $(p_n)_{n\in\IN}$ in $2^\IN$ as follows.
We let $p_n(k)=0$ if $p\langle i,k\rangle=0$ for $i=0,...,n$ and we switch
to $p_n(k)=1$ for all $n\geq m$ if $p\langle m,k\rangle=1$.
Hence the sequence $(p_n)_{n\in\IN}$ converges and the limit contains
a $1$ in those positions $k$ for which there exists $n\in\IN$ with $p\langle n,k\rangle=1$.
Now we compute a non-empty concept class $\CC\in\AA_-$ for this converging sequence $(p_n)_{n\in\IN}$ just as in the proof of Theorem~\ref{thm:VCdim-negative} 
and we set $h(p):=\CC$. Then we obtain
\begin{eqnarray*}
	p\in\GG&\iff&
	(\forall^\infty k)(\exists n)\;p\langle n,k\rangle=1
	\iff(\forall^\infty k)\lim_{n\to\infty}p_n(k)=1\\
    &\iff& \VCdim(h(p))<\infty\iff h(p)\in\PP\LL.
\end{eqnarray*}
This completes the proof.
\end{proof}

We note that $\AA$ is a completely metrizable separable space and the 
Borel structure that we obtain
is the standard Effros Borel structure~\cite{Kec95,BG09}.
The spaces $\AA_-$ and $\AA_+$ are quasi-Polish spaces and
Borel complexity on such spaces can also be defined directly~\cite{dBre13}
(in a way which is equivalent to our approach via realizers).

The effective Borel completeness can also be used to characterize
the Weihrauch complexity of the characteristic function $\chi_{\PP\LL}:\AA\to\{0,1\}$
of $\PP\LL$ that decides PAC learnability. We write $\chi_{\PP\LL}^+$ or
$\chi_{\PP\LL}^-$ if we use $\AA_+$ or $\AA_-$ instead of $\AA$, respectively.

\begin{corollary}[PAC Learnability]
$\chi_{\PP\LL}\equivSW\chi_{\PP\LL}^+\equivSW\LPO'$ and $\chi_{\PP\LL}^-\equivSW\LPO''$.
\end{corollary}
\begin{proof}
From the proof of Theorem~\ref{thm:Borel} we can extract the reductions
$\chi_\FF\leqSW\chi_{\PP\LL}\leqSW\LPO'$ (likewise for $\chi_{\PP\LL}^+$) and $\chi_\GG\leqSW\chi_{\PP\LL}^-\leqSW\LPO''$. The reductions $\LPO'\leqSW\chi_\FF$ and
$\LPO''\leqSW\chi_\GG$ follow as $\LPO^{(n)}=\LPO\circ\lim^{(n-1)}$ is the characteristic function of 
an effective $\Pi^{0}_{n+1}$--set for $n\geq1$ (as $\lim$ is effectively $\Sigma^0_2$--measurable, e.g., by \cite[Proposition~9.1]{Bra05}).
\end{proof}
 
\section{Relative PAC Learning}
\label{sec:RPAC}

In this section we prove that relative PAC learning is computably equivalent
to Weak K\H{o}nig's Lemma, provided that the concept class $\CC\in\AA_+$
is given by positive information, the hypothesis class $\HH\in\AA_-$
is given by negative information and an upper bound of the VC dimension
of $\CC$ is provided as an extra input information.
In order to prove the desired results on relative PAC learning
we need some preparations.

\begin{lemma}[Set of probability distributions]
	\label{lem:prob-distrib}
	The set $\PP$ of probability distributions\linebreak
	${\DD:\IN\times\{0,1\}\to\IR}$ is a c.e.\ closed subset of $\IR^{\IN\times\{0,1\}}$.
\end{lemma}
\begin{proof}
	For simplicity we use the computable encoding
	$f:\IN\to\IN\times\{0,1\}$
	with $f(2n)=(n,0)$ and $f(2n+1)=(n,1)$ that allows us to identify $\IN\times\{0,1\}$
	with $\IN$. Now given a finite sequence of non-empty rational intervals $(a_i,b_i)\In\IR$ with $a_i,b_i\in\IQ$
	for $i=0,...,k$, we need to recognize that the basic open set $(\bigtimes_{i=0}^k(a_i,b_i))\times\IR^\IN$
	contains a probability distribution $\DD:\IN\to\IR$, if this is the case.
	However, this is the case
	if and only if
	\[\sum_{i=0}^ka_i<1\mbox{ and }(\forall i=0,...,k)\;b_i\geq0.\]
	This is even a decidable property.
	It is also easy to see that the set of probability distributions is a closed
	subset, hence it is c.e.\ closed altogether.
\end{proof}

\begin{lemma}[Computability of the risk]
	\label{lem:computability-risk}
	The following functions are computable:
	\begin{enumerate}
		\item $L:\PP\times2^\IN\to\IR,(\DD,h)\mapsto L_\DD(h)$.
		\item $\inf L:\In\PP\times\AA_+\to\IR_>,(\DD,\CC)\mapsto\inf_{h\in\CC}L_\DD(h)$.
	\end{enumerate}
	Here $\inf L$ is defined only for non-empty $C$.
\end{lemma}
\begin{proof}
	If $\DD$ is a probability distribution, then we can determine for every
	value $n\in\IN$ a $k\in\IN$ such that $\sum_{x\leq k,y\in\{0,1\}}\DD(x,y)>1-2^{-n}$.
	This allows us to compute
	\[L_\DD(h)=\sum_{(x,y)\not\in\graph(h)} \DD(x,y)\]
	for every $h\in2^\IN$ with any prescribed precision $2^{-n}$, which proves (1).
	For (2) we note that by Proposition~\ref{prop:seq-positive}
	a non-empty $\CC\in\AA_+$ can be given by a sequence
	$(h_n)_{n\in\IN}$ such that $\CC=\overline{\{h_n:n\in\IN\}}$.
	Then we obtain
	\[\inf L(\DD,\CC)=\inf_{n\in\IN}L_{\DD}(h_n),\]
	which can be computed as a number in $\IR_>$ by (1).
\end{proof}

Next we prove that relative PAC learning can be reduced to Weak K\H{o}nig's Lemma.
We exploit the fact that $\WKL\equivSW\C_\LL$ and we show that
given $\CC$, $\HH$ and an upper bound on the VC dimension of $\CC$, we
can recognize a sufficiently large class of non-learners.

\begin{proposition}
	$\RPAC\leqSW\WKL$.
\end{proposition}
\begin{proof}
	Let $\CC\in\AA_+$, $\HH\in\AA_-$ and $d\in\IN$
	with $\VCdim(\CC)\leq d$ be given.
	Then we can compute a sample complexity for $\CC$
	of the form
	\[m_d\langle i,j\rangle=c\cdot2^i(d\cdot i+j)\]
	for some suitable universal constant $c\in\IN$.
	Hence, by the Fundamental Theorem~\ref{thm:fundamental} we know that there is a learner $A:\SS\to2^\IN$ with $\range(A)\In\HH$
	(even one with $\range(A)\In\CC$, but we do not want to use this condition)
	that PAC learns $\CC$ with the sample complexity $m_d$, i.e., such that
	for every $i,j\in\IN$ and $n\geq m_d\langle i,j\rangle$ and every $\DD\in\PP$
	the learning condition (\ref{eq:PAC}) is satisfied.

	In order to prove that we can find such a learner using $\C_\LL\equivSW\WKL$ and given $\CC,\HH$ and $d$,
	it suffices to prove that we can computably recognize that a function $A:\SS\to2^\IN$
	is not such a learner. But this is equivalent to
	\begin{enumerate}
		\item $(\exists S\in\SS)\;A(S)\not\in\HH$ or
		\item $(\exists \DD\in\PP)(\exists i,j\in\IN)(\exists n\geq m_d\langle i,j\rangle)\;\IP_{S\sim\DD^n}\left[L_\DD(A(S))>\inf_{h\in\CC}L_\DD(h)+2^{-i}\right]> 2^{-j}$.
	\end{enumerate}
	The first condition is clearly c.e.\ open relative to $A$ and $\HH$ as $\HH\in\AA_-$ is
	given with respect to negative information.
	For the second condition, we note that by Lemma~\ref{lem:prob-distrib} the
	set $\PP$ of probability distributions is c.e.\ closed. Hence the existential quantification over $\DD$
	leads to a c.e.\ open condition, if the condition after the quantifiers is c.e.\ open.
	To this end, we note that by Lemma~\ref{lem:computability-risk} we
	can compute $L_\DD(A(S))\in\IR$, given $A$, $\DD$ and $S$, and we can compute $\inf_{h\in\CC}L_{\DD}(h)\in\IR_>$, given $\CC$ and $\DD$.
	Hence, the condition
	\[U=\left\{S\in\SS:L_\DD(A(S))>\inf_{h\in\CC}L_\DD(h)+2^{-i}\right\}\]
	is c.e.\ in its parameters $\DD,\CC$ and $A$. This implies
	that we can compute the probability
	\[\Prob_{S\sim\DD^n}[S\in U]=\sum_{((x_1,y_1),...,(x_n,y_n))\in U}\prod_{i=1}^n\DD(x_i,y_i)\in\IR_<\]
	uniformly in $n\in\IN$,
	since we can sum up all the positive instances that we can confirm.
	This suffices to make the condition $\Prob_{S\sim\DD^n}[S\in U]>2^{-j}$
	c.e.\ open in all its parameters.
	Hence, with the help of $\C_\LL\equivSW\WKL$ we can
	find a learner $A$ that PAC learns $\CC$ relative to $\HH$ for the sample complexity $m_d$.
	This allows us to produce the pair $(A,m_d)$ as output information.
\end{proof}

For the inverse direction of the proof we exploit the fact
that $\WKL\equivSW\widehat{\C_2}$ and for any instance
$(B_k)_{k\in\IN}$ of $\widehat{\C_2}$ we compute
a concept class $\CC\in\AA_+$ of $\VCdim(\CC)\leq1$ together
with a hypothesis class $\HH\in\AA_-$ such that $\CC\In\HH$
and such that any learner $A:\SS\to\HH$ that PAC learns
$\CC$ relative to $\HH$ allows us to reconstruct answers $i\in\C_2(B_k)$
for every $k\in\IN$.

\begin{proposition}
	\label{prop:WKL-IPAC}
	$\WKL\leqSW\RPAC$.
\end{proposition}
\begin{proof}
	Given a sequence $B=(B_k)_{k\in\IN}$ of sets $B_k\in\AA_-(\{0,1\})$ with $B_k\not=\varnothing$
	for all $k\in\IN$ we compute $\CC_B\in\AA_+$ as follows:
	\begin{enumerate}
		\item We add $\widehat{0}$ to $\CC_B$.
		\item We add $h_{k,0}:=0^{k+1}1\widehat{0}$ to $\CC_B$ if $1\not\in B_k$.
		\item We add $h_{k,1}:=0^{k+1}11\widehat{0}$ to $\CC_B$ if $0\not\in B_k$.
	\end{enumerate}
	That is
	\[\CC_B=\{\widehat{0}\}\cup\{h_{k,b}:1-b\not\in B_k\}.\]
	And we compute $\HH_B\in\AA_-$ with $\CC_B\In\HH_B$ as follows:
	\begin{enumerate}
		\item We start with $\HH_B$ that contains everything that potentially goes into $\CC_B$.
		\item We remove $h_{k,1}$ from $\HH_B$ if $1\not\in B_k$.
		\item We remove $h_{k,0}$ from $\HH_B$ if $0\not\in B_k$.
	\end{enumerate}
	That is
	\[\HH_B=\{\widehat{0}\}\cup\{h_{k,b}:k\in\IN,b\in\{0,1\}\}\setminus\{h_{k,b}:b\not\in B_k\}.\]
	Since $B_k\not=\varnothing$ for all $k\in\IN$, it is guaranteed that $\CC_B\In\HH_B$, as
	$\HH_B$ is constructed starting from a set that contains everything that potentially goes into $\CC_B$ by removing only paths that are definitely not added to $\CC_B$.
	We note that the VC dimension of $\CC_B$ is less or equal to $1$.
	Now we assume that $(A,m)\in\RPAC(\CC_B,\HH_B,1)$. In order to determine
	$\C_2(B_k)$ we need to find $b\in B_k$.
	The learner $A$ has to obey the learning condition (\ref{eq:PAC}) in particular
	for the Dirac distribution $\DD$ with $\DD(k+1,1)=1$.
	If $B_k=\{b\}$ with $b\in\{0,1\}$, then $L_\DD(h)=0$ exactly for $h=h_{k,b}\in\CC_B$ and $L_{\DD}(h)=1$ for all other $h\in\HH_B\setminus\{h_{k,b}\}$.
	Hence, for a sample $S=((k+1,1),(k+1,1),...,(k+1,1))\in\SS$ of size $n\geq m\langle 2,2\rangle$
	the learner necessarily has to return $A(S)=h_{k,b}$ with $h_{k,b}(k+2)=b\in B_k$
	in order to satisfy
	\begin{eqnarray}
	\IP_{S\sim\DD^n}\left[L_\DD(A(S))\leq\inf_{h\in\CC_B}L_\DD(h)+\frac{1}{4}\right]> \frac{3}{4}.
	\label{eqn:PAC22}
	\end{eqnarray}
	If $B_k=\{0,1\}$, then we also have $A(S)=h$ for some $h\in\HH_B$ and in any case
	$h(k+2)\in B_k$.
	This proves the reduction $\WKL\equivSW\widehat{\C_2}\leqSW\RPAC$.
\end{proof}

If we modify $\RPAC$ to $\RPAC^{\rm fix}:\In\AA_+\times\AA_-\times\IN\mto\LL$ with
\[\RPAC^{\rm fix}(\CC,\HH,d):=\{A:A\mbox{ PAC learns }\CC\mbox{ relatively to }\HH\mbox{ with }m_d\}\]
with the same domain as $\RPAC$, i.e., if the learner has to work always with the fixed
sample complexity $m_d$,
then all the above proofs still go through. Hence, we have proved the following.
The proofs also work for any fixed finite VC dimension $\geq1$.

\begin{theorem}[Relative PAC learning]
	\label{thm:relative}
	$\RPAC\equivSW\RPACd{d}\equivSW\RPAC^\mathrm{fix}\equivSW\WKL$ with $1<d<\infty$.
\end{theorem}

From Theorem~\ref{thm:relative}
we can draw various positive and negative consequences.
For instance, from the existence of a Kleene tree (a computable binary tree without computable
paths) we can conclude that there are c.e.\ closed concept classes that are included
in co-c.e.\ closed hypothesis class, but do not admit computable learners with computable
sample complexity. However, in Corollary~\ref{cor:PA-learners} 
we are going to formulate a stronger negative result.

On the other hand, the Low Basis Theorem has the consequence expressed in the following corollary (see~\cite{BBP12}).
We recall that a learner $A:\SS\to\HH$ seen as a point in the Cantor space $\LL$
is called {\em low}, if its Turing jump is Turing reducible to the halting problem.
This means that the learner is almost computable in the sense that its oracle power
regarding jumps is not stronger than that of a computable learner.

\begin{corollary}[Lowness of relative PAC learners]
	\label{cor:lowness}
	For every c.e.\ closed concept class $\CC$ of finite VC dimension $d$
	and every co-c.e.\ closed superclass $\HH$ of $\CC$ there is a low learner $A:\SS\to\HH$ that PAC learns $\CC$ relative to $\HH$ with computable sample complexity $m_d$.
\end{corollary}

Finally, it is known~\cite[Theorem~7.2]{BBP12}
that $\WKL$ is complete for non-deterministically computable
problems in the sense of~\cite{Zie07a}. Hence we obtain the following.

\begin{corollary}[Non-deterministic relative PAC learning]
	The relative PAC learning problem $\RPAC$ is complete for the class of non-deterministically
	computable problems.
\end{corollary}

We close this section with a study of the question how the complexity of relative PAC learning
is affected if we do not provide an upper bound on the VC dimension as extra input information,
but if we just keep finiteness of the VC dimension as a promise.
That is, we want to study $\RPACd{\infty}$ in the following.
We prove that it is actually equivalent to Euclidean choice $\C_\IR$.

\begin{theorem}[Relative PAC learning without dimension]
	\label{thm:relative-pac-without-dim}
	$\RPACd{\infty}\equivSW\C_\IR$.
\end{theorem}
\begin{proof}
	By Corollary~\ref{cor:VC} the problem $\VCdim^+_{<\infty}$
	of computing the VC dimension of a class $\CC\In2^\IN$ of finite VC dimension is equivalent to $\C_\IN$.
	Hence, we obtain with Theorem~\ref{thm:relative} and Propositions~\ref{prop:WKL} and \ref{prop:C_R}
	\[\RPACd{\infty}\leqSW\RPAC * \VCdim^+_{<\infty}\leqSW \WKL*\C_\IN\equivSW\C_\IR.\]
	By Proposition~\ref{prop:C_R}
	it suffices to prove $\B\times\widehat{\C_2}\leqSW\RPACd{\infty}$ in order to complete the proof.
	We use a similar construction as in the proof of Proposition~\ref{prop:WKL-IPAC}.
	Given $p\in\dom(\B)$
	and $B=(B_k)_{k\in\IN}\in\dom(\widehat{\C_2})$ we compute
	sets $\CC=\CC_B\cup\CC_p\in\AA_+$ and $\HH=\HH_B\cup\HH_p\in\AA_-$, where
	$\CC_B$ and $\HH_B$ are exactly the sets from the proof of Proposition~\ref{prop:WKL-IPAC}.
	And with $d_p:=\max\{p(n):n\in\IN\}$ we define (with $\cdot$ denoting the {\em concatenation}):
	\begin{itemize}
		\item $\CC_p:=1\cdot\{0,1\}^{d_p}\cdot\widehat{0}$,
		\item $\HH_p:=1\cdot 2^\IN$.
	\end{itemize}
	We note that $d_p=\VCdim(\CC_p)$.
	Given $(A,m)\in\RPACd{\infty}(\CC,\HH)$, we can determine the values of $\C_2(B_k)$ for all $k\in\IN$
	exactly as in the proof of Proposition~\ref{prop:WKL-IPAC}.
	Additionally, with $d=\VCdim(\CC)$ we obtain by the Fundamental Theorem~\ref{thm:fundamental}
	that every $n\in\IN$ with
	\[n\geq\frac{m\langle 0,0\rangle}{C}\geq\frac{M_{d}\langle 0,0\rangle}{C}=d=\VCdim(\CC_B\cup\CC_p)\geq\VCdim(\CC_p)=d_p,\]
	satisfies $n\in\B(p)$. Since we can assume without loss of generality that the constant $C\in\IR$ is computable,
	we can determine a corresponding $n$, given $m$. This concludes the proof of the reduction.
\end{proof}

\section{Proper PAC Learning}
\label{sec:PPAC}

In this section we study proper PAC learning.
If we regard the proper PAC learning setting ($\CC=\HH$) as a special case of the relative one, then the hypothesis class that we consider is described by full information. However, it is also interesting to consider the case where the hypothesis class is only described by positive information.
We also know that if a hypothesis class is properly PAC learnable, then any empirical risk minimizer for $\CC$ properly PAC learns $\CC$. It is therefore interesting to also define an $\ERM$ problem.

\begin{definition}[Empirical risk minimizer]
	We define $\ERM: \AA\setminus\{\varnothing\} \mto\LL$ by
	\[\ERM(\CC):=\{A\in\LL:A\mbox{ is an empirical risk minimizer for }\CC\}.\]
	Once again, we define $\ERM^+$, $\ERM^-$ by replacing $\AA$ with $\AA_+$ or $\AA_-$, respectively.
\end{definition}

We first prove that with full information on the concept class proper
PAC learning is computable, no matter whether an upper bound of 
the VC dimension is explicitly given or just promised to have some specific
finite value.

\begin{theorem}[Proper PAC learning with full information and dimension]
	\label{thm:proper-computable}
	We obtain that
	$\ERM\equivSW\PPAC\equivSW\PPACd{d}\equivSW\id$ are computable for $1<d<\infty$.
\end{theorem}

\begin{proof}
	We claim that we obtain the reduction chain
	\begin{eqnarray}
	\id\leqSW\PPACd{2}\leqSW\PPACd{d}\leqSW\PPAC\leqSW\ERM\leqSW\id
	\label{eqn:PPAC}
	\end{eqnarray}
	for all $1<d<\infty$. The second and the third reduction hold obviously. 
	The first one can be obtained using the computable function $2^\IN\to\AA,p\mapsto\{p\}$ and using the fact that $\id\equivSW\id_{2^\IN}$.
	We need to prove the other reductions.
	
	We start showing that $\ERM$ is computable, which implies the last reduction. 
	Since we have the full information about the hypothesis class $\CC\in\AA$, by Proposition~\ref{prop:trees} 
	we can for any sample $S=((x_1,y_1),...,(x_n,y_n)) \in \SS$ read that information until we know exactly which prefixes of length $x_{\max}+1$ are extendable, where 
	$x_{\max}:=\max\{x_1,...,x_n\}$.
	We can then easily compute the empirical risk for these extendable prefixes, and extend the best of them to a hypothesis $h \in \CC$ using the positive information. This produces an empirical risk minimizer $A$, which solves the problem $\ERM$.

    We now prove $\PPAC\leqSW\ERM$.
	The computed empirical risk minimizer properly PAC learns $\CC$ with $m_{\VCdim(\CC)}$ according to Theorem \ref{thm:fundamental}. We have an upper bound $d$ on $\VCdim(\CC)$ as in input information, so we can simply return $(A,m_d)$ to solve the problem $\PPAC$. This already proves $\PPAC\leqW\ERM$.
	In order to obtain a strong reduction, we need to feed the upper bound
	on the hypothesis class through the problem $\ERM$. However, this can be
	achieved using the computable function \[\AA\times\IN\to\AA,(\CC,d)\mapsto\CC_d:=0^d\cdot1\cdot\CC.\]
		We note that $\VCdim(\CC_d)=\VCdim(\CC)$.
	And if $A_d\in\ERM(\CC_d)$ then we obtain some learner $A\in\ERM(\CC)$,
	which is implicitly defined by 
	\[0^d\cdot1\cdot A((x_1,y_1),...,(x_n,y_n))=A_d((x_1+d+1,y_1),...,(x_n+d+1,y_n)).\]
	The value $d$ itself can be recovered from any value in $\range(A_d)$.
\end{proof}

Now we consider the case of a concept class with full information and only
with the promise that the VC dimension of that class if finite.

\begin{theorem}[Proper PAC learning with full information and without dimension]
	We have $\PPACd{\infty}\equivW \C_\IN$.
\end{theorem}

\begin{proof}
	We have $\PPACd{\infty}\leqW \C_\IN$ because we can return $(A,m_d)$ where $A$ is computed like in the previous proof and $d=\VCdim(\CC)$ is computed using the positive information of $\CC$ with $\C_\IN$.

	For the other direction, we show that $\B \leqW \PPACd{\infty}$, which is sufficient as $\B \equivW \C_\IN$. Given $p\in\dom(\B)$, we write the full information of a hypothesis class $\CC\in\AA$ such that $\VCdim(\CC)=\max(p)$:
	\[
		\CC := 0^{n_1} \cdot \{0,1\} \cdot 0^{n_2} \cdot \{0,1\} \cdot ... \cdot 0^{n_{\max (p)}} \cdot \{0,1\} \cdot \widehat{0},
	\]
	where $n_1,n_2,...,n_{\max(p)}\in\IN$ are chosen as large as needed. In other words, we choose to increase the VC dimension when we encounter higher values in $p$. This is feasible since at any point during the construction of $\CC$, the set is only determined up to a certain depth, so we can allow both values $0$ and $1$ at a later time.

	Then, we obtain an output $(A,m)$ of $\PPACd{\infty}(\CC)$. Like in the proof of Theorem \ref{thm:relative-pac-without-dim}, let $n\in\IN$ be such that $n \geq \frac{m\langle 0,0\rangle}{C}$. Then $n\geq\VCdim(\CC)$, so $n\in\B(p)$.
\end{proof}

We note that the statement cannot be strengthened to a strong equivalence,
as $\PPACd{\infty}$ has inputs for which it can only produce non-computable
outputs, whereas $\C_\IN$ always produces computable outputs. 
However, using the computable map 
\[\AA\times2^\IN\to\AA,(\CC,p)\mapsto\{\langle q,p\rangle:q\in\CC\}\]
together with the construction from the previous proof
we obtain the following corollary, which shows that the strong equivalence
class of $\PPACd{\infty}$ is exactly the cylindrification of $\C_\IN$,
a degree that also contains the limit $\lim_\Delta$ of $\IN^\IN$
with respect to the discrete topology~\cite{BBP12}.

\begin{corollary}
	\label{cor:proper-dim}
	$\PPACd{\infty}\equivSW\id\times\C_\IN\equivSW\lim_\Delta$.
\end{corollary}

Now we study proper PAC learning the only positive or negative information
on the concept class, which turns out to be equivalent to the limit map.
The same holds for the jumps of proper PAC learning with full information.
The later case can also be seen as providing an even weaker limit
information on the concept class.

\begin{theorem}[Proper PAC learning with weaker information]
	\label{thm:PAC-weaker}
	We obtain
	\begin{enumerate}
		\item $\ERM^+\equivSW\PPACd{d}^+\equivSW\PPAC^+\equivSW\lim$ for $1<d\leq\infty$.
		\item $\ERM^-\equivSW\PPACd{d}^-\equivSW\PPAC^-\equivSW\lim$ for $1<d<\infty$.
		\item $\ERM'\equivSW\PPACd{d}'\equivSW\PPAC'\equivSW\lim$ for $1<d<\infty$.
	\end{enumerate}
\end{theorem}
\begin{proof}
(3) is a direct consequence of Theorem~\ref{thm:proper-computable} as $\id'\equivSW\lim$.
(1) and (2) for $1<d<\infty$
 can be proved with a similar reduction chain as in (\ref{eqn:PPAC}), namely
we have
\begin{eqnarray*}
	\lim\leqSW\PPACd{2}^\pm\leqSW\PPACd{d}^\pm\leqSW\PPAC^\pm\leqSW\ERM^\pm\leqSW\lim,
\end{eqnarray*}
where ``$\pm$'' is to be read such that we can either have a ``$+$'' everywhere or a ``$-$'' everywhere. Then the second and third reductions are obvious again.
The last mentioned reduction holds as $\ERM^\pm\leqSW\ERM'\equivSW\lim$ by (3)
and Proposition~\ref{prop:convert}. The reduction $\PPAC^\pm\leqSW\ERM^\pm$
can be proved as in the proof of Theorem~\ref{thm:proper-computable}, noting
that the function 
\[\AA_\pm\times\IN\to\AA_\pm,(\CC,d)\mapsto\CC_d:=0^d\cdot1\cdot\CC\]
used there is computable also for negative or positive information alone.
We still need to prove $\lim\leqSW\PPACd{2}^\pm$ and we use
$\lim\equivSW\widehat{\LPO}$ for this purpose, which holds by Proposition~\ref{prop:lim}.
For a given input $p=\langle p_0,p_1,...\rangle$ of $\widehat{\LPO}$ we can compute
the concept classes
\begin{itemize}
\item  $\CC^+=\{\widehat{0}\}\cup\{0^k1\widehat{0}: (\exists i\in\IN)\ p_{k}(i)=1\}\in\AA_+$,
\item $\CC^-=\{\widehat{0}\}\cup\{0^k1\widehat{0}:(\forall i\in\IN)\ p_{k}(i)=0\}\in\AA_-$.
\end{itemize}
The positive information on $\CC^+$ can be computed by adding the
hypothesis $h_k=0^k1\widehat{0}$ whenever a $1$ is encountered in $p_k$
and the negative information on $\CC^-$ can be computed by removing
$h_k$ from $\{\widehat{0}\}\cup\{0^k1\widehat{0}:k\in\IN\}$ whenever
a $1$ is encountered in $p_k$.
We note that we have $\VCdim(\CC^\pm)=1$. 

Now we obtain a learner $A$ and a sample complexity $m$ for $\CC^\pm$ with $\PPACd{2}^\pm(\CC^\pm)$. 
Using $m$, we compute a sample size $n\geq m\langle 2,2\rangle$ such that the accuracy and the confidence errors are both $\frac{1}{4}$. 
Then for a sample $S=((k,1),...,(k,1))$ of size $n$ the learner $A$ has to satisfy 
the learning condition as in (\ref{eqn:PAC22}).
Thus, by considering the Dirac distribution $\DD(k,1)=1$ we obtain
that $h=A(S)$ satisfies $h(k)=1$ if and only if there is a one in $p_k$ in the case of $\CC^+$ and if and only if there is no one in $p_k$ in the case of $\CC^-$.
In both cases we can recover the information on $\widehat{\LPO}(p)$ from $h=A(S)$,
which proves $\widehat{\LPO}\leqSW\PPACd{2}^\pm$. 

It remains to prove (1) for the case $d=\infty$. But now we can proceed
as in the proof of Theorem~\ref{thm:relative-pac-without-dim} and use
$\C_\IN$ to compute the VC dimension (see Corollary~\ref{cor:VC}). We obtain with Proposition~\ref{prop:lim}
\[\PPACd{2}^+\leqSW\PPACd{\infty}^+\leqSW\PPAC^+*\VCdim_{<\infty}^+\leqSW\lim*\C_\IN\equivSW\lim,\]
which completes the proof using (1) for the case $d=2$.
\end{proof}

We note that the reduction of $\lim$ to the problem of proper PAC learning
only requires to evaluate the sample complexity $m$ on a single input $m\langle 2,2\rangle$.
Hence, the entire non-computability is included in the learner and we obtain the
following non-uniform result that strengthens the known result that there are RER concept
classes of VC dimension $1$ that do not admit computable learners, proved
by Agarwal et al.~\cite[Theorem~9]{AAB+20}.

\begin{corollary}[Proper learners of the degree of the halting problem]
	\label{cor:proper-halting-problem}
	There is a c.e.\ closed (and also another co-c.e.\ closed) concept class $\CC\In2^\IN$ of
	VC dimension $1$ such that every learner $A:\SS\to\CC$ that properly learns $\CC$
	is necessarily of the Turing degree of the halting problem. 
\end{corollary}
	
By Theorem~\ref{thm:proper-computable} such a concept class $\CC$ cannot be 
computably closed (i.e., simultaneously
c.e.\ closed and co-c.e.\ closed). And the Turing degree of the halting problem
provides the maximal complexity that can be achieved. 
We call a closed concept class {\em limit computable}, if it is computable in the space $\AA'$.
This set of classes includes both, c.e.\ closed and co-c.e.\ closed concept classes.
	
\begin{corollary}[Limit computable risk minimizers]
	\label{cor:limit-computable-risk-minimizer}
	Every limit computable closed concept class $\CC\In2^\IN$ of finite VC dimension
	admits a limit computable empirical risk minimizer $A:\SS\to\CC$.
\end{corollary}

The following lemma tells us that PAC learning problems are idempotent, i.e.,
we can use one instance to solve two instances. 

\begin{lemma}[Idempotency of PAC learning]
	\label{lem:product}
	Given two non-empty concept classes $\CC_i\In2^\IN$ of VC dimension
	$\VCdim(\CC_i)=d_i$ for $i\in\{0,1\}$, the concept class $\CC=\langle \CC_0,\CC_1\rangle$
	satisfies $\VCdim(\CC)=d_0+d_1$. Given a learner $A:\SS\to2^\IN$ for $\CC$ with sample
	complexity $m$, we can compute learners $A_i:\SS\to2^\IN$ for $\CC_i$ with the same sample complexity $m$. Additionally, if $A:\SS\to\CC$ is proper for $\CC$ then the learners $A_i:\SS\to\CC_i$ are proper for $\CC_i$ for $i\in\{0,1\}$.
\end{lemma}
\begin{proof}
	The class $\CC=\langle\CC_0,\CC_1\rangle=\{\langle p,q\rangle:p\in\CC_0\mbox{ and }q\in\CC_1\}$ is defined such that the hypotheses from $\CC_0$ fill the even positions
	and those from $\CC_1$ fill the odd positions.
	Let $B\In \IN$ be a finite set shattered by $\CC$. Then $\CC_i$ shatters $B_i:=\{n\in\IN: 2n+i\in B\}$ for $i\in\{0,1\}$ and $|B|=|B_0|+|B_1|$, proving $\VCdim(\CC)\leq d_0+d_1$. Conversely, if $\CC_i$ shatters the finite set $B_i\In\IN$ for $i\in\{0,1\}$, then $\CC$ shatters the union $B$ of the disjoint sets $2 B_0$ and $2B_1 +1$
	with $|B|=|B_0|+|B_1|$, showing the reverse inequality. This proves $\VCdim(\CC)=d_0+d_1$.
	
	Let now $A:\SS\to2^\IN$ be a learner for $\CC$ with sample complexity $m$.
	Then we define learners $A_i:\SS\to2^\IN$ by
	\[A_i((x_1,y_1),...,(x_n,y_n)):=\pr_{i+1} A((2x_1+i,y_1),...,(2x_n+i,y_n))\]
	for all samples $S=((x_1,y_1),...,(x_n,y_n))\in\SS$, where $\pr_1:\IN^\IN\to\IN^\IN,\langle p,q\rangle\mapsto p$ denotes the projection
	on the even positions and $\pr_2:\IN^\IN\to\IN^\IN,\langle p,q\rangle\mapsto q$
	the projection on the odd positions. Since the learner $A$ has to satisfy the
	learning condition (\ref{eq:PAC}) for $\CC$ and
	all distributions $\DD$, it does so, in particular
	for those distributions that put all the weight only on even or only on odd positions.
	But this implies that $A_i$ satisfies the learning condition for $\CC_i$
	and all distributions with the same sample complexity $m$. In order to see this,
	we note that if $\DD_i$ is some distribution, then there is a corresponding
	distribution $\DD$  that puts the same weights as $\DD_i$ but only on
	even or odd positions depending on whether $i=0$ or $i=1$ and we obtain $\inf_{h\in\CC_i} L_{\DD_i}(h)=\inf_{h\in\CC} L_\DD(h)$.
	It is also clear that $\LL\to\LL\times\LL,A\mapsto(A_0,A_1)$ is computable.
	
	It is clear that the learner $A_i$ defined above satisfies $\range(A_i)\In\CC_i$,
	if $\range(A)\In\CC$.
\end{proof}

We note that in (2) and (3) of Theorem~\ref{thm:PAC-weaker}
the case $d=\infty$ cannot be added,
as we obtain a different complexity in this case.

\begin{proposition}
	\label{prop:PPAC-ininity}
	We obtain $\PPACd{\infty}^-\equivSW\PPACd{\infty}'\equivSW\lim\times\C_\IN'$.
\end{proposition}
\begin{proof}
Since $\B\equivW\C_\IN$ we obtain  $\id\times\B\equivSW\id\times\C_\IN$ and
hence $\lim\times\B'\equivSW\lim\times\C_\IN'$. Thus it suffices to prove
\[\lim\times\B'\leqSW\PPACd{\infty}^-\leqSW\PPACd{\infty}'\equivSW\lim\times\C'_\IN.\]
The second reduction and the equivalence hold by Proposition~\ref{prop:lim} and Corollary~\ref{cor:proper-dim}.
For the other reduction we use Corollary~\ref{cor:VCdim-negative}, 
which implies $\B'\leqSW\PPACd{\infty}^-$, as we can compute an upper
bound on the VC dimension using the sample complexity $m$ and the constant $C$
by $n\geq\frac{m\langle 0,0\rangle}{C}\geq\VCdim(\CC)$
as in the proof of Theorem~\ref{thm:relative-pac-without-dim}.
On the other hand, Theorem~\ref{thm:PAC-weaker} implies
$\lim\leqSW\PPAC^-\leqSW\PPACd{\infty}^-$. 
The reduction $\lim\times \B'\leqSW\PPACd{\infty}^-$ then follows
as Lemma~\ref{lem:product} implies $\PPACd{\infty}^-\times\PPACd{\infty}^-\leqSW\PPACd{\infty}^-$
because $\AA_-\times\AA_-\to\AA_-,(\CC_0,\CC_1)\mapsto\langle \CC_0,\CC_1\rangle$
is computable.
\end{proof}

The construction in the proof of Theorem~\ref{thm:PAC-weaker}
also proves that converting positive into negative information and vice versa
is equivalent to the limit operation, even restricted to closed sets of VC dimension one.
This is because it is sufficient to transfer $\CC^+$ from that proof into $\CC^-$ or
vice versa to compute $\widehat{\LPO}$. The other direction in the following
corollary follows from Proposition~\ref{prop:convert}.

\begin{corollary}
	\label{cor:conversion}
	The identities $\id:\AA_+\to\AA_-$ and $\id:\AA_-\to\AA_+$ are strongly
	Weihrauch equivalent to $\lim$, even restricted to closed sets of VC dimension one.
\end{corollary}

\section{Improper PAC Learning}
\label{sec:IPAC}

In this section we want to analyze the problem of improper PAC learning
and it variants. Since improper PAC learning is just PAC learning relative to $2^\IN$,
we directly obtain 
\[\IPAC^+\leqSW\RPAC\equivSW\WKL\]
with the help of Theorem~\ref{thm:relative}.
But is the upper bound $\WKL$ tight and is $\IPAC^+$ equivalent to it?
We will see that this is not the case, but it turns out that $\IPAC^+$
is in the finitary DNC range, a statement that we will make precise
and that has a number of interesting consequences. 

For our study it is helpful to analyze the proof of Theorem~\ref{thm:Fundamental-computable}, the computable
version of the Fundamental Theorem by Delle Rose et al.\ in order
to understand improper PAC learning.
Hence, our first goal is to extract the uniform content of their proof.

Firstly, we need to recall the No Free Lunch Theorem~\cite[Theorem~5.1]{SB14},
which is one of the basic theorems
in machine learning and it roughly states that every learner fails badly on some task.
Agarwal, Ananthakrishnan, Ben-David, Lechner, and Urner derived a computable
version of this theorem~\cite[Lemma~19]{AAB+20} which allows one to compute
for every learner and sample size $n$ a specific uniform distribution for some finite set 
of size $k\geq 2n$ on which the learner badly fails. However, the computability
assumption on the learner $A$ in their formulation of the theorem is superfluous 
and their proof actually shows the following formally stronger uniform result.

\begin{theorem}[No Free Lunch Theorem]
\label{thm:NFL}
Given an arbitrary learner $A:\SS\to2^\IN$ and $n\in\IN$ as well as natural numbers
$x_1<...<x_k$ for some $k\geq 2n$, we can computably find a 
function $g:\{x_1,...,x_k\}\to\{0,1\}$ such that
\[\IP_{S\sim\DD^n}\left[L_\DD(A(S))\geq\frac{1}{8}\right]\geq\frac{1}{7}\]
for the uniform distribution $\DD$ over $\{(x_1,g(x_1)),...,(x_k,g(x_k))\}$.
\end{theorem}

We can now more or less exactly follow the
proofs given by Delle Rose et al.~\cite[Theorem~11]{DRKRS23}
and by Sterkenburg~\cite[Lemma~9]{Ste22a}
in order to obtain the following result that captures the uniform content 
of the proof of Theorem~\ref{thm:Fundamental-computable}.
For completeness we go through sketches of their proofs in order to see 
that everything holds uniformly.
By $\WW_k$ we denote the set of functions $f:\IN^{k}\to\{0,1\}^k$
and by $\WW:=\bigsqcup_{k\in\IN}\WW_k$ the {\em space of witness functions},
where $\bigsqcup$ stands for the {\em coproduct}.
This space is represented such that a name $\langle k,p\rangle$ comes
with the explicit value $k$ and a description $p\in\IN^\IN$ of $f$.

\begin{theorem}[Witnesses and improper learners]
\label{thm:witness-learner}
We obtain the following:
\begin{enumerate}
	\item Given a witness $f\in\WW$ of some non-empty concept class $\CC$ (that is not given), we
	can computably find a learner $A:\SS\to2^\IN$ and a sample complexity $m\in\IN^\IN$
	such that $A$ improperly learns $\CC$ with $m$. 
\item Given a learner $A:\SS\to2^\IN$ and a sample complexity $m\in\IN^\IN$
such that $A$ improperly learns some non-empty concept class $\CC$ (that is not given)
with sample complexity $m$, we can computably find a witness $f\in\WW$ for $\CC$.
\end{enumerate}
\end{theorem}
\begin{proof}
\begin{enumerate}
	\item We sketch the proof given by Delle Rose et al.~\cite[Theorem~11]{DRKRS23} in order to see that it yields the desired uniformly computable solution. Given the witness $f\in\WW$ of some non-empty concept class $\CC\In2^\IN$ (which is not given)
	we have $k\in\IN$ and the witness function $f:\IN^{k+1}\to\{0,1\}^{k+1}$.
	In particular, we know $\VCdim(\CC)\leq k$. Following Delle Rose et al.\ we call
	a hypothesis $h:\IN\to\{0,1\}$ {\em good for $f$} 
	if $\supp(h):=h^{-1}\{1\}$ is finite
	and for all $x_1<...<x_{k+1}<\max\supp(h)$ we obtain
	\[f(x_1,...,x_{k+1})\not=(h(x_1),...,h(x_{k+1})).\]
	By $\HH_f$ we denote the class of all hypotheses which are good for $f$.
	The goal is to compute a pair $(A,m)$ such that $A:\SS\to2^\IN$ is a learner that
	learns $\CC$ relative to $\HH:=\CC\cup\HH_f$ with sample complexity $m:\IN\to\IN$.
	As proved by Delle Rose et al.\ we obtain $\VCdim(\HH)\leq k+1$. And given
	$f$, we can compute an empirical risk minimizer $A:\SS\to2^\IN$ for $\HH_f$
	since given $S=((x_1,y_1),...,(x_n,y_n))\in\SS$ with $M=\max\{x_1,...,x_n\}$
	we can just go through all finitely many good sequences $h$ with $\supp(h)\In\{0,...,M\}$ in order to find one with minimal error $L_S(h)$.
	As proved by Delle Rose et al.\ this $h$ is also an empirical risk minimizer
	for $\HH$ and hence according to Theorem~\ref{thm:fundamental} a learner for $\CC\In\HH$ with sample complexity $m_{k+1}$. Hence, we can return $(A,m_{k+1})$
	as a result upon input of $f\in\WW$.
    \item We now sketch the proof given by Sterkenburg~\cite[Lemma~9]{Ste22a}
    in order to see that it yields the desired uniformly computable solution.
    Given a learner $A:\SS\to2^\IN$ in $\LL$ that learns a non-empty concept
    class $\CC\In2^\IN$ (that is not given) with sample complexity $m:\IN\to\IN$.
    we have for $n=m\langle 4,3\rangle$ that
    $\IP_{S\sim\DD^n}\left[L_\DD(A(S))\leq\inf_{h\in\CC}L_\DD(h)+\frac{1}{16}\right]>1- \frac{1}{8}$ and hence, in particular,
    \begin{eqnarray}
    	\IP_{S\sim\DD^n}\left[L_\DD(A(S))\geq\inf_{h\in\CC}L_\DD(h)+\frac{1}{8}\right]<\frac{1}{7}
    \end{eqnarray}
    for every probability distribution $\DD$.
	For $k=2n$ and all natural numbers $x_1<...<x_k$ we can then find by the 
	No Free Lunch Theorem~\ref{thm:NFL} a function $g:\{x_1,...,x_k\}\to\{0,1\}$ such that
	for the uniform distribution $\DD$ over $\{(x_1,g(x_1)),...,(x_k,g(x_k))\}$ we have
	\[\IP_{S\sim\DD^n}\left[L_\DD(A(S))\geq\frac{1}{8}\right]\geq\frac{1}{7},\]
	which implies that $\graph(g)\not\In\graph(h)$ for every $h\in\CC$. In particular,
	we obtain a witness function $f:\IN^k\to\{0,1\}^k$ for $\CC$ if we define 
	\[f(x_1,...,x_k):=(g(x_1),...,g(x_k)).\]	
    This completes the proof.\qedhere
\end{enumerate}	
\end{proof}

We note that if we start with some pair $(A,m)$ and
apply first (2) and then (1), then the proof yields a new pair $(A',m')$ where
$m'$ is of the form $m'=m_{k+1}$ for some $k\in\IN$. Hence, $m'$ is always computable.

\begin{corollary}[Computable sample complexities]
	\label{cor:computable-sample}
	Given some learner $A:\SS\to2^\IN$ and some $m:\IN^\IN\to\IN^\IN$ such that $A$ learns some non-empty concept class	$\CC\In2^\IN$ (that is not given) with sample 
	complexity $m$, we can computably find another learner $A':\SS\to2^\IN$  and 
	a computable sample complexity $m':\IN\to\IN$ such that $A'$ learns
	$\CC$ with sample complexity $m'$.
\end{corollary}

In this sense, we can always obtain some computable sample complexity
with a corresponding learner. In passing, we note that the uniform formulation
of the results by Delle Rose et al.\ and Sterkenburg given in Theorem~\ref{thm:witness-learner}
does not only imply Theorem~\ref{thm:Fundamental-computable}, but also 
the following relativized version, which also generalizes the qualitative part of the
classical Fundamental Theorem of Statistical Learning, i.e., Theorem~\ref{thm:fundamental}.

\begin{corollary}[Relativized Fundamental Theorem of Statistical Learning]
For a non-empty concept class $\CC\In2^\IN$ and a Turing degree ${\mathbf a}$
the following are equivalent: $\CC$ admits an improper learner $A:\SS\to2^\IN$
that is computable from ${\mathbf a}$ (with computable sample complexity $m$),
if and only if $\CC$ admits a witness function that is computable from ${\mathbf a}$.
\end{corollary}

Theorem~\ref{thm:witness-learner} tells us 
that computing a witness for a concept class given with some 
upper bound on its VC dimension is as hard as computing an improper learner 
of the same class with corresponding sample complexity. 
In order to express this more formally, we introduce the problem
of computing witnesses.

\begin{definition}[Witness problem]
	We consider the {\em witness problem} $\WIT:\In \AA \times\IN\mto\WW$, defined by
	\[\WIT(\CC,d):=\{f:\IN^{k+1}\to\{0,1\}^{k+1}:f\mbox{ is a witness function for $\CC$}\}\]
	and $\dom(\WIT)=\{(\CC,d):\VCdim(\CC)\leq d$ and $\CC\neq\varnothing\}$.
	We use analogous notations with the upper index $+$ or $-$ if we replace
	the space $\AA$ by either $\AA_+$ or $\AA_-$, respectively.
\end{definition}

We immediately obtain the following corollary of Theorem~\ref{thm:witness-learner}.

\begin{corollary}[Improper learning and witnesses]
\label{cor:witness-learner}
We obtain the following equivalences:
$\IPAC\equivSW\WIT$, $\IPAC^+\equivSW\WIT^+$, and $\IPAC^-\equivSW\WIT^-$.
\end{corollary}

Studying the witness problem is in many cases easier than studying the
improper PAC learning problem itself. For instance, we directly obtain the following.

\begin{proposition}[Computing witnesses]
	\label{prop:witness-computable}
	$\WIT^-$ is computable.
\end{proposition}
\begin{proof}
Given a concept class $\CC\in\AA_-$ and $d\geq\VCdim(\CC)$,
we know that for $k=d+1$ there exists a $(k-1)$--witness function 
$f:\IN^{k}\to\{0,1\}^k$ for $\CC$. 
Given natural numbers $x_1<...<x_k$ we consider $M:=\max\{x_1,...,x_k\}$
and $W_b:=\{w\in\{0,1\}^{M+1}:(w(x_1),...,w(x_k))\not=b\}$ for every $b\in\{0,1\}^k$.
Since $\CC\in\AA_-$ is a compact set, given $(x_1,...,x_k)$ and $b$, the property
\[\CC\In\bigcup_{w\in W_b}w2^\IN\]
is c.e.\ in all parameters by~\cite[Theorems~4.9 and 4.10]{BP03}. And the property has to be satisfied for some of the finitely many $b\in\{0,1\}^k$ since there exists
a $(k-1)$--witness function. Hence we can find some $b\in\{0,1\}^k$ for which
it holds and we can set $f(x_1,...,x_k):=b$ to define a witness function.
This yields an algorithm to compute $f$, given $\CC\in\AA_-$ and $d$.
\end{proof}

The following corollary is immediate with the help of Corollaries~\ref{cor:VC} and \ref{cor:VCdim-negative}.

\begin{corollary}
	$\IPAC$ and $\IPAC^-$ are computable, $\IPACd{\infty}\equivW\C_\IN$ and $\B'\leqSW\IPACd{\infty}^-\leqW\C_\IN'$.
\end{corollary}

The algorithm in the proof of Proposition~\ref{prop:witness-computable}
also yields that we can always find a witness function for the exact
VC dimension, if that is given. Hence, we obtain the following.

\begin{corollary}[Improper PAC learnability of co-c.e.\ closed classes]
Every non-empty co-c.e.\ closed concept class $\CC$ satisfies $\VCdim(\CC)=\eVCdim(\CC)$.
In particular, every non-empty co-c.e.\ closed concept class $\CC$ 
has a computable improper learner (with computable sample complexity) if and only
if it is of finite VC dimension.
\end{corollary}

\begin{figure}[htb]
	\begin{center}
		\begin{tikzpicture}[scale=0.9,auto=left,thick,every node/.style={fill=blue!20},
			PAC/.style ={fill=violet!20}]
			\node (WKL) at (-4,7) {$\DNC_2\equivSW\WKL$};
			\node (DNC3) at (-4,5) {$\DNC_3$};
			\node (DNC4) at (-4,4) {$\DNC_4$};
			\node[style={fill=none}] (DNCD) at (-4,3) {$\vdots$};
			\node (ACC3) at (-2,4) {$\ACC_3$};
            \node (ACC4) at (-2,3) {$\ACC_4$};			
			\node[style={fill=none}] (ACCD) at (-2,0.5) {$\vdots$};
			\node (DNCS) at (-4,2) {$\DNC_*$};
			\node (WWKL) at (0.5,6) {$\WWKL$};
			\node (PA) at (-6,-0.5) {$\PA$};
			\node (DNC) at (-4,0.5) {$\DNC_\IN$};
			\node (ACCN) at (-2,-0.5) {$\ACC_\IN$};			
			\node (LLPO) at (-2,5) {$\LLPO$};
			\node[PAC] (IPAC) at (-7.5,5) {$\IPAC^+\equivSW\WIT^+$};
			\node[PAC] (IPACd) at (-6,4) {$\IPACd{d}^+$};
			\node[PAC] (IPAC2) at (-6,3) {$\IPACd{2}^+$};	
			\node[PAC] (IPACPS) at (-7.5,2) {$\IPAC^+_{-\mathrm{sample}}$};	
			\draw[->] (WKL) edge (WWKL);
			\draw[->] (WKL) edge (DNC3);
			\draw[->] (DNC3) edge (DNC4);
			\draw[->] (DNC4) edge (DNCD);
			\draw[->] (DNCD) edge (DNCS);
			\draw[->] (LLPO) edge (ACC3);
			\draw[->] (ACC3) edge (ACC4);
			\draw[->] (DNC3) edge (ACC3);
			\draw[->] (DNC4) edge (ACC4);
			\draw[->] (WWKL) edge (LLPO);
			\draw[->] (DNCS) edge (PA);
			\draw[->] (DNCS) edge (DNC);
			\draw[->] (WWKL) edge  [bend angle=30, bend left] (DNC);
			\draw[->] (WWKL) edge (LLPO);
			\draw[->] (WKL) edge (IPAC);
			\draw[->] (IPAC) edge (IPACd);
			\draw[->] (IPACd) edge (IPAC2);
			\draw[->] (IPAC2) edge (DNCS);
			\draw[->] (IPAC) edge (IPACPS);
			\draw[->] (IPACPS) edge (PA);
			\draw[->] (DNC) edge (ACCN);			
			\draw[->] (ACC4) edge (ACCD);
			\draw[->] (ACCD) edge (ACCN);						
		\end{tikzpicture}
		\caption{The vicinity of $\DNC$ in the Weihrauch lattice.}
		\label{fig:DNC}
	\end{center}
\end{figure}
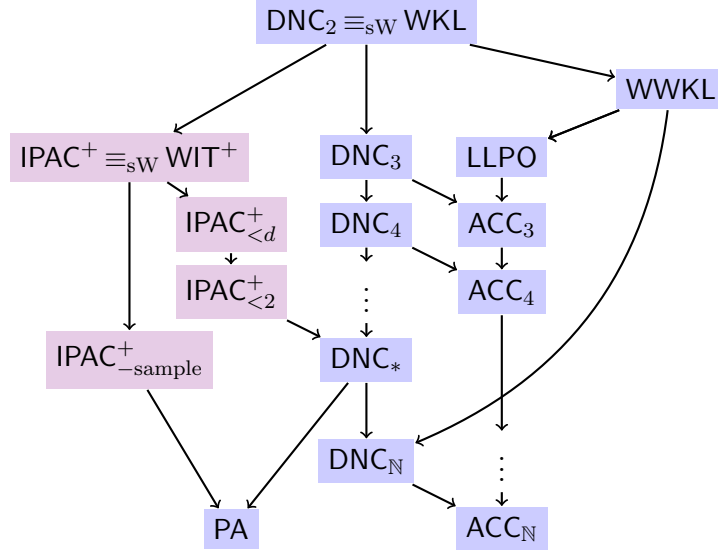

Now we want to study improper PAC learning from positive information $\IPAC^+$.
In order to give the reader some orientation, we provide in Figure~\ref{fig:DNC} a 
finer resolution of the relevant part of the Weihrauch lattice from Figure~\ref{fig:Benchmark}. 
Some more details can be found in \cite{BHK17a}.

We have already mentioned that $\IPAC^+\leqW\WKL$ follows from
Theorem~\ref{thm:relative}. Another independent proof can be obtained
from the following observation.

\begin{proposition}[Computing the set of witnesses]
\label{prop:set-of-witness-functions}
The map 
\[F:\AA_+\to\AA_-(\WW), \CC\mapsto\{f\in\WW:f\mbox{ is a witness function for }\CC\}\] 
is well-defined and computable. In particular, the set of witnesses for a 
given concept class is a closed subset of $\WW$.
\end{proposition}
\begin{proof}
One immediately obtains this result as for $\CC\in\AA_+$ and $f:\IN^k\to\{0,1\}^k\in\WW$
\[f\not\in F(\CC)\iff(\exists h\in\CC)(\exists x_1<...<x_k)\;f(x_1,...,x_k)=(h(x_1),...,h(x_k)).\]
Existential quantification over sets $\CC\in\AA_+$ with c.e.\ open predicates is c.e.\ open
relative to $\CC$, which yields the desired result.
\end{proof}

We note that the set of witness functions might be empty. But if the
VC dimension of $\CC$ is given as an input information, then we know some $k$
such that $\WW_k$ definitely contains a witness function for $\CC$ and
hence $\IPAC^+\equivSW\WIT^+\leqSW\WKL$ follows (we implicitly use here
that $\WW_k$ is computably isomorphic to $2^\IN$ and hence
$\CC_{\WW_k}\equivSW\CC_{2^\IN}\equivSW\WKL$).

Next we want to separate $\IPAC^+$ from $\WKL$ by looking at witness functions.
In fact, since $\WKL\equivSW\widehat{\LLPO}$ it suffices to prove $\LLPO\nleqW\WIT^+$ for this purpose.
We prove a result which is stronger in two respects.
For one, we even consider $\ACC_n$ for general $n\geq 2$ (note that $\ACC_2\equivSW\LLPO$) and secondly we provide a topological separation.
Here $\leq_{\mathrm W}^*$ denotes the topological version of Weihrauch reducibility,
where the functions $H$ and $K$ just have to be continuous.

\begin{proposition}
	$\ACC_n\not\leq_{\mathrm W}^*\WIT^+$ for all $n\geq2$.
\end{proposition}
\begin{proof}
Let us assume for a contradiction that $\ACC_n\leq_{\mathrm W}^*\WIT^+$ for some $n\geq2$.
Then there are continuous $H,K$ such that $H\langle\id,GK\rangle$ is a realizer
for $\ACC_n$ whenever $G$ is a realizer for $\WIT^+$.
We now consider a name $p\in\IN^\IN$ of the set $A=\{0,...,n-1\}\in\AA_-(\{0,...,n-1\})$.
This set is a legitimate input for $\ACC_n$ and hence $K(p)$ has to
be a name for some concept class $\CC\in\AA_+$ together with
some upper bound $d\geq\VCdim(\CC)$ of the VC dimension of $\CC$.
A realizer for $\WIT^+$ can produce a name of any witness $f\in\WW$ 
for $\CC$. We only consider the case that $f\in\WW_k$ for $k=2^{n}(d+1)$.
The set $W\In\WW_k$ of possible witnesses for the set $\CC$ is closed, as a finite
prefix of a witness suffices to confirm that it fails, see Proposition~\ref{prop:set-of-witness-functions}.
Since $W$ is a closed subset of the compact space $(\{0,1\}^k)^{\IN^k}$
it is itself compact.
We can assume that the representation of $\WW$ is such that the functions
are represent bijectively via some natural coding. Hence, we can identify $W$
with the set of its names.
Thus, $H$ is uniformly continuous on the compact set $\{p\}\times W$
and some finite prefix $w\prefix p$ suffices for $H$ to produce a natural
number output in $\{0,...,n-1\}$ for every $f\in W$.
Additionally, we can assume by continuity of $K$ that $w$ is long enough such that $K$
already produces the upper bound $d$ of the VC dimension on the output side.
For every $i\in\{0,...,n-1\}$ there is an extension $p_i\in\IN^\IN$ of $w$
that is a name of the set $A_i:=\{0,...,n-1\}\setminus\{i\}$.
By choice of $w$ these $p_i$ are mapped to names $K(p_i)$ of 
concept classes $\CC_i\in\AA_+$ with $\VCdim(\CC_i)\leq d$.
Hence, the set $\CC':=\CC\cup\bigcup_{i=0}^{n-1}\CC_i$ is a union of $n+1$ sets
all of VC dimension $\leq d$ and thus
$\VCdim(\CC')\leq 2^{n}(d+1)-1$ by a well-known formula
for the VC dimension of unions.\footnote{For concept classes $\HH_i$ with $\VCdim(\HH_i)\leq d$ one obtains $\VCdim(\HH_1\cup\HH_2)\leq 2d+1$, which inductively yields $\VCdim(\bigcup_{i=1}^{n+1}\HH_i)\leq2^n(d+1)-1$, even though much better
estimations are known~\cite[Excercise~6.11]{SB14}.}
Hence, there is a witness function $f\in\WW_k$ for $\CC'$ as $k=2^{n}(d+1)$.
Since $\CC\In\CC'$, this witness function also satisfies $f\in W$ and hence
$H$ will map $p$ together with a name for $f$ to some number $i\in\{0,...,n-1\}$.
Since $\CC_i\In\CC$, the function $f$ is also a witness function for $\CC_i$
and hence a possible output of a realizer of $\WIT^+$ upon input of a name of $\CC_i$ and $d$.
By choice of $w$ the same number $i$ will also be the result of $H$ upon
input of $p_i$ and $f$.	However, this result is incorrect as $i\not\in A_i$.
Hence, the reduction fails and we have a contradiction!
\end{proof}

The argument in the proof cannot be generalized to the countable case as
$\ACC_\IN\leqW\WIT^+$. In fact, we get a stronger result by proving
that even $\DNC_*$ is reducible to $\WIT^+$.

\begin{proposition}
	\label{prop:DNC-WIT}
	$\DNC_*\leqSW\WIT^+$.
\end{proposition}
\begin{proof}
Given a sequence $A=((A_{n,m})_{n\in\IN})_{m\geq2}$ of sets $A_{n,m}\in\AA_-(m)$ 
as an input
to $\DNC_*$, we have that at most one $i\in\{0,...,m-1\}=m$ is missing in each $A_{n,m}$.
The goal is to find some pair $(m,p)\in\IN\times\IN^\IN$ 
such that $(\forall n)\;p(n)\in A_{n,m}$. 
This goal has to be achieved with the help of $\WIT^+(\CC,d)$
for some suitable concept class $\CC\in\AA_+$ together with some bound $d\geq\VCdim(\CC)$ that need to be computed from the input data.
In fact, we will only attempt to find pairs of the form $(2^k,p)$ with $(\forall n)\;p(n)\in A_{n,{2^k}}$ 
for sets with indexes that are powers of $2$. This is because
the witness problem yields functions of type $f:\IN^k\to\{0,1\}^k$.
There are exactly $2^k$ possible outputs of such functions and we denote
by $\bin_k(i)\in\{0,1\}^k$ the $(i+1)$--st such binary word of length $k$,
numbered using the binary representation of $i\in\{0,...,2^k-1\}$ of length $k$.
In fact, we can compute
\[\CC=\{\widehat{0}\}\cup\{0^{2^{\langle n,k\rangle}}\cdot\bin_k(i)\cdot\widehat{0}:i\not\in A_{n,{2^k}},i<2^k,k\geq1,n\in\IN\}\in\AA_+\]
from the sequence $A$ of sets. 
We use the standard definition of Cantor pairs that satisfies $\langle n,k\rangle\geq k$.
Hence, the prefixes $0^{2^{\langle n,k\rangle}}$ 
together with the fact that there is at most one $i<2^k$ missing in each set $A_{n,2^k}$ ensure
that two blocks of the form $\bin_k(i)$ in two hypotheses in $\CC$ for different 
values $\langle n,k\rangle$ never overlap and hence $\VCdim(\CC)\leq1$. 
Given $\CC$ and the upper bound $1$
of the VC dimension as an input to $\WIT^+$, 
we obtain some witness function $f\in\WIT^+(\CC,1)$
of type $f:\IN^k\to\{0,1\}^k$ for some $k\geq1$.
Now for every $n\in\IN$ we determine some value $p(n)\in\{0,...,2^k-1\}$ such that
\[f(2^{\langle n,k\rangle},2^{\langle n,k\rangle}+1,...,2^{\langle n,k\rangle}+k-1)=\bin_k(p(n))\]
and since $f$ is a witness function for $\CC$ we can conclude that 
$(\forall n)\;p(n)\in A_{n,2^k}$.
Hence $\langle 2^k,p\rangle$ is a suitable solution for $\DNC_*(A)$, which proves the reduction.	
\end{proof}

Since the reduction yields a concept class of VC dimension $1$, we obtain
the following immediate corollary.

\begin{corollary}
	$\DNC_*\leqSW\IPACd{2}^+\leqSW\IPACd{d}^+\leqSW\IPAC^+\lW\WKL$ for all $0<d<\infty$.
\end{corollary}

As $\DNC_2\equivSW\WKL$, we have in particular 
\begin{eqnarray}
	\DNC_*\leqSW\IPAC^+\lW\DNC_2.
\label{eqn:finitary-DNC-range}
\end{eqnarray}
This is what we mean by saying that $\IPAC^+$ is {\em in the finitary
DNC range} and it implies that in several respects $\IPAC^+$ behaves
like a problem that is similar (but not equivalent) to $\DNC_n$ for some $n\in\IN$. 

Since $\PA$ and $\DNC_\IN$ are two (in fact incomparable) lower bounds of $\DNC_*$,
we obtain a number of important consequences. The following follows
from Proposition~\ref{prop:finitary-DNC}.

\begin{corollary}
	$\PA\lW\IPACd{2}^+\nleqW\WWKL^{(n)}$ for all $n\in\IN$.
\end{corollary}

The proof of Proposition~\ref{prop:finitary-DNC} implicitly uses a classical result of 
Jokusch and Soare~\cite[Corollary~5.4]{JS72}, which shows that the set of 
PA degrees has measure 0.
In some sense this rules out any possibility of a meaningful probabilistic
algorithm for improper PAC learning from positively represented concept classes.
In particular, improper PAC learning for such classes is not a probabilistic problem and not
Las Vegas computable in the sense of~\cite{BGH15a}.

\begin{corollary}
	Improper PAC learning $\IPAC^+$ from positive information is not Las Vegas computable (not even restricted to concept classes of VC dimension $1$).
\end{corollary}

By \cite[Proposition~6.6]{BHK17a} it is known that $\WKL\leqW\C'_\IN*\PA$.
Hence, $\C_\IN'$ gives us also
an upper bound on how far $\IPAC^+$ is away from $\WKL$.

\begin{corollary}
	\label{cor:WKL-CNP-IPAC}
	$\WKL\leqSW\C_\IN'*\IPAC^+$.
\end{corollary}

We can also get some conclusion regarding the problem $\IPACd{\infty}^+$.
Similarly as for the proof of Proposition~\ref{prop:PPAC-ininity} we can 
use Lemma~\ref{lem:product} to obtain
an interval that encloses $\IPACd{\infty}^+$.

\begin{corollary}
	$\DNC_*\times\C_\IN\leqW\IPAC^+\times\C_\IN\leqW\IPACd{\infty}^+\leqW\IPAC^+*\C_\IN\leqW\C_\IR$.
\end{corollary}

Finally, we briefly want to discuss the third question (C) that we have raised in the 
introduction, namely what happens if we omit the sample complexity on the output
side of our PAC learning problems. We only consider $\IPAC^+$ as this is the
weakest natural non-computable PAC learning problem that we have
considered. 
Let us denote by $\IPAC^+_{-\mathrm{sample}}$ the version of $\IPAC^+$
where we remove te sample complexity from the output.
How much weaker is $\IPAC^+_{-\mathrm{sample}}$ than $\IPAC^+$?
If the output is just a learner without any guaranteed sample complexity,
than the learner could do anything on samples of small size and hence
the set of outputs $\IPAC^+_{-\mathrm{sample}}(\CC,d)$ for any 
arbitrary input $(\CC,d)$ in the domain is a dense subset of $\LL$.
Hence $\IPAC^+_{-\mathrm{sample}}(\CC,d)$ is densely realized
in the sense of \cite{BHK17a} and by \cite[Proposition~4.3]{BHK17a}
we obtain the following conclusion.

\begin{corollary}
	\label{cor:ACCN-IPAC}
$\ACC_\IN\nleqW\IPAC^+_{-\mathrm{sample}}$.
\end{corollary} 

This applies analogously to all other PAC learning problems that we have considered,
if we remove the sample complexity on the output side.
In this sense they become very weak, as in a well-defined sense
no discrete information can be extracted from them.
This corresponds to the intuition that from a practical perspective
a learner without sample complexity is useless.

On the positive side, we can extract some information from the proof
of Theorem~\ref{thm:witness-learner}. When the witness function
is computed from the pair $(A,m)$ in that proof, the sample complexity
$m$ is only evaluated on a single input in order to determine $n=m\langle 4,3\rangle$. 
If the concept class $\CC\in\AA_+$ is also given, then we can just guess the number $n$
instead and verify the witness function $f:\IN^k\to\{0,1\}^k$ using
Proposition~\ref{prop:set-of-witness-functions}. If the witness function
fails, we replace $n$ by $n+1$ and start over again.
This yields the following result.

\begin{proposition}
Given a non-empty concept class $\CC\in\AA_+$ of finite VC dimension together
with some learner $A:\SS\to2^\IN$ that improperly learns $\CC$,
we can determine a witness function $f\in\WW$ for $\CC$ in a finite mind
change computation. Hence, we can also find 
in a finite mind change computation another learner $A':\SS\to2^\IN$
together with a (computable) sample complexity $m$ such that $A'$ improperly
learns $\CC$ with $m$.
\end{proposition}

This yields the following corollary.

\begin{corollary}
	$\IPAC^+\leqW\C_\IN*\IPAC^+_{-\mathrm{sample}}$.
\end{corollary}

With Corollary~\ref{cor:WKL-CNP-IPAC}, the fact that $\C_\IN'\leqW\C_\IN*\C_\IN'$
(which follows for instance from \cite[Theorem~7.2]{SV23})
and the fact that $\PA\leqW \C_\IN*f\TO\PA\leqW f$ (which holds
as finite mind change computations can only produce solutions that are 
Turing reducible to the respective instances)
we obtain the following.

\begin{corollary}
	$\PA\leqW\IPAC^+_{-\mathrm{sample}}$ and
	$\WKL\leqW\C_\IN'*\IPAC^+_{-\mathrm{sample}}$.
\end{corollary}

In some sense this result is another way of confirming the fact that the
essential contribution to the non-computability of PAC learning
comes from the learner, not from the sample complexity. 
In fact, all our results would remain the same if we would replace the sample
complexity $m\in\IN^\IN$ by the natural number $m\langle 4,3\rangle\in\IN$,
as the sample complexity was never used on any input larger than $\langle 4,3\rangle$
in any of our proofs (and since we can assume $m\langle 4,3\rangle\geq m\langle 2,2\rangle$, where $\langle 2,2\rangle$ was the only other input on which $m$ was evaluated).

We can also formulate a corollary that addresses the complexity of learners
that belong to c.e.\ closed concept classes.
Our result strengthens the result that there are RER concept classes $\CC$ of VC dimension $1$
that do not admit computable learners, which was proved by Sterkenburg~\cite[Theorem~3]{Ste22a}.
This following result is a corollary of Proposition~\ref{prop:DNC-WIT}
and the proof of Theorem~\ref{thm:witness-learner}. Since we 
only need to exploit the sample complexity on a single input $m\langle 4,3\rangle$
within that proof the entire non-computability is included in the learner
and once again we do not need to mention the sample complexity here.

\begin{corollary}[Improper learners of PA degree]
	\label{cor:PA-learners}
	There is a c.e.\ closed concept class $\CC\In2^\IN$ of VC dimension $1$
	such that every witness function for $\CC$
	and every learner $A:\SS\to2^\IN$ that improperly learns $\CC$ 
	is of PA degree.
	This class $\CC$ is necessarily of infinite effective VC dimension.
\end{corollary}

A positive result in the other direction is given in
Corollary~\ref{cor:lowness}, which also applies to improper learners
(which are learners relative to $\HH=2^\IN$).

\section{Conclusions}
\label{sec:conclusion}

In this article we have studied computability aspects of the problem: {\em given} a concept class with upper bound on its VC dimension, {\em find} a learner with a sample complexity for it. We emphasize
that we have looked at this question from a uniform perspective, which means that
we focus on the complexity of finding the learner (with a corresponding sample complexity).
The previously existing non-uniform results rather focus on the mere existence of computable
learners and conditions that can guarantee such an existence. 
Some of our results suggest that we can look at 
witness functions as additional input information on the concept classes.
They contain more information than a mere upper bound of the VC dimension.
In the diagram in Figure~\ref{fig:costs} we summarize our results from this perspective.

\begin{figure}[htb]
	\usetikzlibrary {arrows.meta,automata,positioning,shadows}
	\begin{tikzpicture}
		[shorten >=1pt,node distance=3cm,on grid,>={Stealth[round]},thick,{font=\small},bend angle=10,
		state/.style={circle, draw, minimum size=2cm},
		full/.style ={draw=blue!70,very thick, fill=blue!30, text=black},
		positive/.style ={draw=green!70,very thick, fill=green!30, text=black},
		negative/.style ={draw=red!70,very thick, fill=red!30,text=black},
		learn/.style ={draw=orange!70,very thick, fill=orange!30,text=black}]	
		
		\node[state,positive]  (AP)                  {$\CC\in\AA_+$};
		\node[state,positive]  (APd) [right=of AP]   {$\normalsize\CC\in\AA_+\atop d>\VCdim(\CC)$};
		\node[state,positive]  (APw) [right=of APd]  {$\CC\in\AA_+\atop\mbox{ \tiny witness }f$};
		\node[state,full]      (A)   [below=of AP]   {$\CC\in\AA$};
		\node[state,full]      (Ad)  [right=of A]    {$\CC\in\AA\atop d>\VCdim(\CC)$};
		\node[state,full]      (Aw)  [right=of Ad]   {$\CC\in\AA\atop\mbox{ \tiny witness }f$};
		\node[state,negative]  (AN)  [below=of A]    {$\CC\in\AA_-$};
		\node[state,negative]  (ANd) [right=of AN]   {$\CC\in\AA_-\atop d>\VCdim(\CC)$};
		\node[state,negative]  (ANw) [right=of ANd]  {$\CC\in\AA_-\atop\mbox{ \tiny witness $f$}$};
		\node[state,learn]     (ERM) [left=of A]     {$\CC\in\AA\atop A\in\ERM(\CC)$};
		\node[state,learn]     (IPAC) [right=of Aw]  {$(A,m)\in\atop\IPAC(\CC,d)$};
		\node[state,learn]     (PPAC) [below=of ANd]  {$(A,m)\in\atop\PPAC(\CC,d)$};
		
		\path[->] (AP) edge [bend right] node [left] {$\lim$} (A)	
		(A) edge [bend right] node [left] {$\id$} (AN)
		(AN) edge [bend right] node [right] {$\lim$} (A)
		(A) edge [bend right] node [right] {$\id$} (AP)
		(APd) edge [bend right] node [left] {$\lim$} (Ad)	
		(Ad) edge [bend right] node [left] {$\id$} (ANd)
		(ANd) edge [bend right] node [right] {$\lim$} (Ad)
		(Ad) edge [bend right] node [right] {$\id$} (APd)
		(APw) edge [bend right] node [left] {$\lim$} (Aw)	
		(Aw) edge [bend right] node [left] {$\id$} (ANw)
		(ANw) edge [bend right] node [right] {$\lim$} (Aw)
		(Aw) edge [bend right] node [right] {$\id$} (APw)	          
		(AP) edge [bend left] node [above] {$\C_\IN$} (APd)
		(APd) edge [bend left] node [above] {$\WKL$} (APw)
		(APw) edge [bend left] node [below] {$\id$} (APd)
		(APd) edge [bend left] node [below] {$\id$} (AP) 
		(A) edge [bend left] node [above] {$\C_\IN$} (Ad)
		(Ad) edge [bend left] node [below] {$\id$} (A)               
		(AN) edge [bend left] node [above] {$\C'_\IN$} (ANd)
		(ANd) edge [bend left] node [below] {$\id$} (AN) 
		(A) edge node [above] {$\id$} (ERM)
		(Aw) edge node [above] {$\id$} (IPAC)
		(APw) edge node [above right] {$\id$} (IPAC)
		(ANw) edge node [below right] {$\id$} (IPAC)
		(ERM) edge [bend angle=30, bend right] node [below left] {$\C_\IN$} (PPAC)
		(PPAC) edge [bend angle=30, bend right] node [below right] {$\id$} (IPAC)
		(Ad) edge [bend angle=35, bend left] node [below right,pos=0.8] {$\id$} (PPAC)
		(AN) edge node [left] {$\lim\times\C_\IN'$} (PPAC)
		(ANd) edge node [left] {$\lim$} (PPAC)
		;
		\path[<->] 
		(Ad) edge node [above] {$\id$} (Aw)
		(ANd) edge node [above] {$\id$} (ANw)
		;
	\end{tikzpicture}
	\caption{Costs of PAC learning for closed concept classes $\CC\In2^\IN$ of finite VC dimension.}
	\label{fig:costs}
\end{figure}

The graph shows the computational costs that are sufficient to obtain
an empirical risk minimizer or a proper or improper PAC learner with sample complexity, respectively.
These costs can be calculated with this graph starting from certain
types of information on the concept class by taking the most cost efficient path
in the diagram, where all costs need to be composed with the compositional
product $*$ of the respective Weihrauch degrees. 
According to our results, this typically yields a 
pretty tight upper bound for the respective problem.
The overall assumption is that we are dealing with non-empty closed concept classes
$C\In2^\IN$ of finite VC dimension.

The readers who are interested in reverse mathematics~\cite{Sim99,Hir15,DM22} might ask
how our results translate into classifications in second-order arithmetic.
There is no direct meta theorem that would allow us to translate the results
automatically, even though some translations are possible~\cite{Uft21}. 
Reverse mathematics is a proof-theoretic approach that asks
for axiom systems that are sufficient and perhaps also necessary to prove
certain results in second-order arithmetic. 

\begin{figure}[htb]
	\begin{small}
		\begin{tabular}{cc}
			{\bf Weihrauch degree} & {\bf reverse mathematics system} \\[0.1cm]\hline
			&\\[-0.3cm]
			$\lim$ & $\ACA_0$\\
			$\WKL$ & $\WKL_0^*$ \\
			$\WWKL$ & $\WWKL_0^*$ \\
			$\DNC_\IN$ & $\mathsf{DNR}$\\
			$\C_\IN$ & $\mathbf{I\Sigma^0_1}$ \\
			$\C_\IN'$ & $\mathbf{I\Sigma^0_2}$
		\end{tabular}
	\end{small}
	\caption{Weihrauch complexity versus reverse mathematics.}
	\label{fig:reverse}
\end{figure}

To the best of our knowledge,
the Fundamental Theorem of Statistical Learning has not yet been analyzed
in reverse mathematics. In many cases Weihrauch complexity results turn out
to be resource sensitive and uniform versions of corresponding reverse mathematics
classifications. In these cases the dictionary from the table in Figure~\ref{fig:reverse}
provides the correspondences between Weihrauch degrees and systems
in reverse mathematics.\footnote{The systems $\WKL_0^*$ and $\WWKL_0^*$ refer
to the versions of these axiom systems without $\Sigma^0_1$--induction $\mathbf{I\Sigma^0_1}$.}

A question that we have not answered in our study is whether improper PAC learning
from positive information is dependent on the VC dimension of the concept class.

\begin{question}
Does $\IPACd{d}^+\lW\IPACd{d+1}^+$ hold true for all $d>1$?
\end{question}

Hence, we did not separate all the problems shown in the diagram
of Figure~\ref{fig:DNC}. We did not even prove that $\IPAC^+$
is not reducible to $\DNC_n$ for $n>2$.

\begin{question}
	Are there $d>1$ and $n>2$ such that $\IPACd{d}^+\leqW\DNC_n$ hold true? 
\end{question}
 
While it is clear that $\IPAC^+\leqW\DNC_2$ and $\DNC_n\nleqW\IPAC^+$, we have
no results that could help to answer the above question.

\subsection*{Acknowledgments}
The first author was funded by the National Research Foundation of South Africa (NRF) -- grant number 151597.
The second author was funded by the Erasmus Programme of the European Union.
We also acknowledge support by the Alexander von Humboldt Foundation.
We would like to thank Emmanuel Rauzy for helpful discussions.

%
%
%
\bibliographystyle{abbrv}
\bibliography{C:/Users/Documents/Spaces/Research/Bibliography/lit}

\newpage
\section*{Appendix: A result on effective VC dimension}

In \cite[Question~6.1]{KK25} Kattermann and Krapp have 
asked whether there exists a PAC learnable class $\CC$
with a computable improper learner and a non-computable hypothesis $h\in\CC$.
The purpose of this appendix is to answer this question in the affirmative
(using Theorem~\ref{thm:Fundamental-computable}).

\begin{proposition}
	There exists a computable sequence $(h_i)_{i\in\IN}$ of hypotheses 
	$h_i:\IN\to\{0,1\}$ that converges to a non-computable $h:\IN\to\{0,1\}$
	of the degree of the halting problem and such that $\HH=\{h\}\cup\{h_i:i\in\IN\}$
	forms a hypothesis class of effective VC dimension $1$. 
\end{proposition}
\begin{proof}
	Let $r:\IN\to\IN$ be a computable enumeration of the halting problem $K=\range(r)$.
	Let $R_i:=\{r(0),...,r(i)\}$, $m_i:=\max(R_i)$ and $v_{i,j}\in\{0,1\}^*$ defined by
	\[v_{i,j}:=\left\{
	\begin{array}{ll}
		1\cdot 0 & \mbox{if $j\in R_i$}\\
		0\cdot 1 & \mbox{otherwise}
	\end{array}\right.\]
	for all $i\in\IN$ and $j\leq m_i$. We define a computable sequence $(h_i)_{i\in\IN}$
	of the form 
	\[h_i:=u_i\cdot\widehat{0}\mbox{ with }
	u_i:=v_{i,0}\cdot 0^{k_{i,0}}\cdot v_{i,1}\cdot 0^{k_{i,1}}\cdot...v_{i,m_i-1}\cdot 0^{k_{i,m_{i}-1}}\cdot v_{i,m_i}\]
	with even values $k_{i,j}\in\IN$ that are defined in an inductive construction over $i\in\IN$.
	The induction in stage $i=0$ starts with $k_{0,j}:=0$ for all $j<m_0$.
	In the induction step $i\to i+1$ we suppose that we have already defined $h_0,...,h_i$
	with corresponding words $u_0,...,u_i$ and even values $k_{0,j},...k_{i,j}$. 
	We need to define $h_{i+1}$. In the case that $r(i+1)\in R_i$, we set $h_{i+1}:=h_i$
	and hence $k_{i+1,j}:=k_{i,j}$ for all $j<m_{i+1}=m_i$.
	Now we assume that $r(i+1)\not\in R_i$ and we need to distinguish two subcases.
	In the case that $r(i+1)>m_i$, we set $k_{i+1,j}:=0$ for all $j$ with $m_i\leq j<m_{i+1}$
	and $k_{i+1,j}:=k_{i,j}$ for $j<m_i$.
	The interesting case is the case of $\iota:=r(i+1)< m_i$.
	In this case we set $k_{i+1,j}:=0$ for all $j$ with $\iota\leq j<m_i$, 
	we set $k_{i+1,j}:=k_{i,j}$ for $j<\iota-1$,
	and we select a necessarily even value for $k_{i+1,\iota-1}$ such that 
	\[|u_i|=|v_{i+1,0}\cdot 0^{k_{i+1,0}}\cdot...\cdot v_{i+1,\iota-1}\cdot 0^{k_{i+1,\iota-1}}|.\]
	In other words, $u_{i+1}$ has a common prefix with $u_i$ including the block
	$v_{i+1,\iota-1}$, but $v_{i+1,\iota}\not=v_{i,\iota}$. After the common block $v_{i+1,\iota-1}$ the word $u_{i+1}$ is filled 
	up with zeros up to the length of $u_i$ and then it continues with $v_{i+1,\iota}$.
    Since $K$ is an infinite set, it is guaranteed that the sequence $(|u_i|)_{i\in\IN}$
    is unbounded and grows monotonically. 

	This construction ensures that if bits are changed from $h_i$ to $h_{i+1}$ 
	in a position $n<|u_i|$, then they can only be changed from $1$ to $0$ and never
	from $0$ to $1$. 
	And if a bit in position $n<|u_i|$ is changed from $1$ to $0$, 
	then all other bits of $h_i$ in positions $k$ with $n<k<|u_i|$ are also set to $0$.
	New bits $1$ can only be introduced in positions $n\geq|u_i|$.
	This property allows us to compute a witness
	$w:\IN^2\to\{0,1\}^2$ for $\{h_i:i\in\IN\}$ as follows.
	Given two natural numbers $n<k$, we first determine some $m\in\IN$ with $|u_m|>k$ and then we set
	\[w(n,k):=\left\{\begin{array}{ll}
		(0,1) & \mbox{if $(\exists i\leq m)\;(h_i(n),h_i(k))=(1,1)$}\\
		(1,1) & \mbox{otherwise}
		\end{array}\right..\]
	We explain for both cases why this choice of values for the witness function works:
	\begin{enumerate}
	\item If there exists $i\leq m$ with $(h_i(n),h_i(k))=(1,1)$,
    then the bit $1$ in position $n$ of $h_i$ 
    cannot be replaced for later $j>i$ by $h_j(n)=0$
    without also replacing the bit $1$ in position $k$ by $h_j(k)=0$ according to the inductive
    construction above. Likewise, there can also not be an earlier hypothesis $h_j$ with
    $j<i$ such that $h_j(n)=0$ and $h_j(k)=1$, because this would imply 
    $|u_j|>k$ and hence the $0$ in the position $n$ could not have been changed to $h_i(n)=1$.
    Hence, 
    $(h_i(n),h_i(k))\not=(0,1)$ for all $i\in\IN$.
    \item If $(h_i(n),h_i(k))\not=(1,1)$ for all $i\leq m$, then we claim that this
    holds for all $i\in\IN$. 
    This is because if $h_m(n)=0$ or $h_m(k)=0$, then the corresponding bit 
    cannot be set to $1$ later. Since the pattern $(1,1)$ does also not occur earlier, we can conclude that
    $(h_i(n),h_i(k))\not=(1,1)$ for all $i\in\IN$.
	\end{enumerate}
	Altogether we have described a computation of a witness function $w:\IN^2\to\{0,1\}^2$  for $\{h_i:i\in\IN\}$ and hence for the closure
	$\HH=\{h\}\cup\{h_i:i\in\IN\}$ of the set $\{h_i:i\in\IN\}$.
	The hypothesis $h$ is in the closure, since the construction ensures that the sequence $(h_i)_{i\in\IN}$ converges to $h$, which has the Turing degree of the halting problem, because
	$h$ is of the form
	\[h=v_{0}\cdot 0^{2k_0}\cdot v_1\cdot 0^{2k_1}\cdot v_2\cdot 0^{2k_2}...,\]
	where $v_j=1\cdot 0\iff j\in K$ and $v_j=0\cdot 1$ otherwise. 
\end{proof}

\end{document}